\documentclass{gtart_h}

\def\ifplaintex{\expandafter\ifx\csname documentclass\endcsname\relax}

\def\gtp{{\mathsurround=0pt\it $\cal G\mskip-2mu$eometry \&\ 
$\cal T\!\!$opology $\cal P\!$ublications}}  

\def\recd{{\small Received:\qua\receiveddate\ifx\reviseddate\relax
\else\qquad Revised:\qua\reviseddate\fi\par}} 

\def\lognumber#1{\def\thelognumber{#1}}
\def\volumenumber#1{\def\thevolumenumber{#1}}
\def\volumeyear#1{\def\thevolumeyear{#1}}
\def\papernumber#1{\def\thepapernumber{#1}}
\def\pagenumbers#1#2{\def\startpage{#1}\def\finishpage{#2}}
\def\published#1{\def\publishdate{#1}}

\def\received#1{\def\receiveddate{#1}}

\def\accepted#1{\def\accepteddate{#1}}

\long\def\asciiabstract#1{\long\def\theasciiabstract{#1}}

\let\\\par\let\thelognumber\relax\let\thevolumenumber\relax
\let\thepapernumber\relax\let\thevolumeyear\relax\let\startpage\relax
\let\finishpage\relax\let\publishdate\relax\let\receiveddate\relax
\let\reviseddate\relax\let\accepteddate\relax\let\theasciititle\relax
\let\theasciiauthors\relax
\let\theasciiabstract\relax

\let\theasciiemail\relax

\ifplaintex
\font\logobig=cmssbx10 scaled 3836
\font\logomed=cmssbx10 scaled 2557
\else
\font\logobig=cmssbx10 scaled 4200
\font\logomed=cmssbx10 scaled 2800
\fi

\long\def\makeagttitle{   
\count0=\startpage
\agt\hfill      
\hbox to 45truept{\vbox to 0pt{\vglue -13truept{\logomed A\kern -.37em{\logobig 
T}\kern -.38em G}\vss}\hss}
\break
{\small Volume \thevolumenumber\ (\thevolumeyear)
\startpage--\finishpage\nl
Published: \publishdate}

\vglue .25truein

{\parskip=0pt\leftskip 0pt plus
1fil\def\\{\par\smallskip}{\Large\bf\thetitle}\par\medskip} \vglue
0.05truein

{\parskip=0pt\leftskip 0pt plus 1fil\def\\{\par}{\sc\theauthors}
\par\medskip}%
 
\vglue 0.03truein 

{\small\leftskip 25truept\rightskip 25truept{\bf Abstract}\stdspace\theabstract

{\bf AMS Classification}\stdspace\theprimaryclass
\ifx\thesecondaryclass\relax\else; \thesecondaryclass\fi\par
{\bf Keywords}\stdspace \thekeywords\par}\vglue 7truept

}   

\ifplaintex
\font\phead=cmsl9 scaled 950
\font\pnum=cmbx10 scaled 913
\font\pfoot=cmsl9 scaled 950
\headline{\vbox to 0pt{\vskip -4.5mm\line{\small\phead\ifnum
\count0=\startpage ISSN 1472-2739 (on-line) 1472-2747 (printed)
\hfill {\pnum\folio}\else\ifodd\count0\def\\{ }%
\ifx\theshorttitle\relax\thetitle\else\theshorttitle\fi\hfill{\pnum\folio}
\else\def\\{ and }{\pnum\folio}\hfill\ifx\theshortauthors\relax\theauthors
\else\theshortauthors\fi\fi\fi}\vss}}
\footline{\vbox to 0pt{\vglue 0mm\line{\small\pfoot\ifnum\count0=\startpage
\copyright\ \gtp\hfill\else
\agt, Volume \thevolumenumber\ (\thevolumeyear)\hfill\fi}\vss}}
\else
\font=cmsl9 scaled 1050
\font\lnum=cmbx10 
\font=cmsl9 scaled 1050
\makeatletter
\def\@oddhead{{\small\lhead\ifnum\count0=\startpage ISSN 1472-2739 
(on-line) 1472-2747 (printed)\hfill {\lnum\number\count0}\else\ifodd\count0
\def\\{ }\ifx\theshorttitle\relax \thetitle \else\theshorttitle\fi\hfill
{\lnum\number\count0}\else\def\\{ and }{\lnum\number\count0}
\hfill\ifx\theshortauthors\relax 
\theauthors\else\theshortauthors\fi\fi\fi}}\def\@evenhead{\@oddhead}
\def\@oddfoot{\small\lfoot\ifnum\count0=\startpage\copyright\ \gtp\hfill\else
\agt, Volume \thevolumenumber\ (\thevolumeyear)\hfill\fi}
\def\@evenfoot{\@oddfoot}
\makeatother
\fi
\let\maketitlepage\makeagttitle
\let\makeshorttitle\maketitlepage
\let\maketitle\maketitlepage

\newwrite\gtoutfile
\long\gdef\makeheadfile{  
{\def\\{, }\def\s{ }
\immediate\openout\gtoutfile head.xxx
\immediate\write\gtoutfile{Proxy-for: \ifx\theasciiauthors\relax
\theauthors\else\theasciiauthors\fi\s<\ifx\theasciiemail\relax\theemail\else\theasciiemail\fi>}
\immediate\write\gtoutfile{\noexpand\\}
\immediate\write\gtoutfile{Authors: \ifx\theasciiauthors\relax
\theauthors\else\theasciiauthors\fi}
{\def\\{ }\immediate\write\gtoutfile{Title: \ifx\theasciititle\relax
\thetitle\else\theasciititle\fi}}
\immediate\write\gtoutfile{Subj-class: GT or SG, GR etc}
\immediate\write\gtoutfile{MSC-class: \theprimaryclass\ifx\thesecondaryclass\relax\else, \thesecondaryclass\fi}
\immediate\write\gtoutfile{Journal-ref: Algebr. Geom. Topol. \thevolumenumber\s
(\thevolumeyear) \startpage-\finishpage}
\immediate\write\gtoutfile{Comments: Published by Algebraic and
Geometric Topology at}
\immediate\write\gtoutfile{\s\s\s  http://www.maths.warwick.ac.uk/agt/AGTVol\thevolumenumber/agt-\thevolumenumber-\thepapernumber.abs.html}
\immediate\write\gtoutfile{\noexpand\\}
\immediate\write\gtoutfile{}
\ifx\theasciiabstract\relax
\immediate\write\gtoutfile{\theabstract}\else
\immediate\write\gtoutfile{\theasciiabstract}\fi
\immediate\write\gtoutfile{}
\immediate\write\gtoutfile{\noexpand\\}
\immediate\write\gtoutfile{}
\immediate\closeout\gtoutfile}}  

\def\maketitlepage{\makeagttitle\makeheadfile}
\let\makeshorttitle\maketitlepage
\let\maketitle\maketitlepage

\lognumber{19}
\volumenumber{4}
\volumeyear{2004}
\papernumber{19}
\published{8 June 2004}
\pagenumbers{347}{398}
\received{17 March 2004}
\accepted{26 March 2004}

\usepackage{amsmath,amssymb,amscd,enumerate}
\input rlepsf

\def\sh#1{\goodbreak{\bf #1}\nobreak
\addcontentsline{toc}{subsection}{#1}}

\def\bbz{\mathbb Z}
\def\bbk{\mathbb K}
\def\bbr{\mathbb R}
\def\bbq{\mathbb Q}
\def\bbc{\mathbb C}
\let\x\times
\let\ox\otimes
\let\op\oplus
\def\SA{\mathcal{A}}
\def\SK{\mathcal{K}}
\def\SR{\mathcal{R}}
\def\SC{\mathcal{C}}
\def\SM{\mathcal{M}}
\def\SE{\mathcal{E}}
\def\SB{\mathcal{B}}
\def\BQ{\mathbf{Q}}
\def\BZ{\mathbf{Z}}
\let\d\delta
\let\a\alpha
\let\b\beta
\let\e\epsilon
\let\G\Gamma
\let\Om\Omega
\let\g\gamma
\let\s\sigma
\let\k\kappa
\let\th\theta
\let\lra\longrightarrow
\let\hra\hookrightarrow
\let\thra\twoheadrightarrow
\let\p\partial
\let\<\langle
\let\>\rangle
\let\sbq\subseteq
\let\cd\cdot
\let\wt\widetilde
\let\tl\widetilde
\def\ngth{\negthickspace}
\let\un\underline
\let\ov\overline
\def\supth{{\rm th}}

\newtheorem*{5.4}{Theorem}
\newtheorem*{7.4}{Corollary}
\newtheorem*{9.5}{Theorem 9.5}
\newtheorem{thm}{Theorem}[section]
\newtheorem{lem}[thm]{Lemma}
\newtheorem{prop}[thm]{Proposition}
\newtheorem{cor}[thm]{Corollary}

\theoremstyle{definition}
\newtheorem{defn}[thm]{Definition}
\newtheorem{ex}[thm]{Example}
\newtheorem{rem}[thm]{Remark}
\newtheorem*{rems}{Remarks}
\newtheorem*{ques}{Question}
\newtheorem*{con}{Conjecture}

\newcommand{\gs}[1]{\ensuremath{G^{(#1)}}}
\newcommand{\gns}[1]{\ensuremath{G^{(n#11)}}}
\newcommand{\sk}{\ensuremath{S^3\backslash K}}
\newcommand{\ltri}{\ensuremath{\triangleleft}}
\renewcommand{\l}{\ell}
\newcommand{\rks}{\operatorname{rk_{\SK}}}
\newcommand{\id}{\operatorname{id}}
\newcommand{\sss}{\ensuremath{S^1\ngth\x\ngth S^1\ngth\x\ngth
S^1}}

\newcommand{\arf}{\operatorname{Arf}}
\newcommand{\aut}{\operatorname{Aut}}
\newcommand{\Ext}{\operatorname{Ext}}
\newcommand{\Hom}{\operatorname{Hom}}
\newcommand{\image}{\operatorname{image}}
\newcommand{\Ker}{\operatorname{Ker}}
\newcommand{\rank}{\operatorname{rank}}
\newcommand{\genus}{\operatorname{genus}}
\newcommand{\spin}{\ensuremath{\mathrm{Spin}}}

\begin{document}
\title{Noncommutative knot theory}
\author{Tim D. Cochran}
\address{Department of Mathematics, Rice University\\6100 Main Street, 
Houston, Texas 77005-1892, USA}
\email{cochran@rice.edu}

\begin{abstract}
The classical \emph{abelian} invariants of a knot are the Alexander
module, which is the first homology group of the the unique infinite
cyclic covering space of $S^3-K$, considered as a module over the
(commutative) Laurent polynomial ring, and the Blanchfield linking
pairing defined on this module. From the perspective of the knot
group, $G$, these invariants reflect the structure of
$G^{(1)}/G^{(2)}$ as a module over $G/G^{(1)}$ (here $G^{(n)}$ is the
$n^\supth$ term of the derived series of G). Hence any phenomenon
associated to $G^{(2)}$ is invisible to abelian invariants. This paper
begins the systematic study of invariants associated to solvable
covering spaces of knot exteriors, in particular the study of what we
call the \emph{$n^\supth$ higher-order Alexander module},
$G^{(n+1)}/G^{(n+2)}$, considered as a
$\mathbb{Z}[G/G^{(n+1)}]$--module.  We show that these modules share
almost all of the properties of the classical Alexander module. They
are torsion modules with higher-order Alexander polynomials whose
degrees give lower bounds for the knot genus.  The modules have
presentation matrices derived either from a group presentation or from
a Seifert surface. They admit higher-order linking forms exhibiting
self-duality. There are applications to estimating knot genus and to
detecting fibered, prime and alternating knots. There are also
surprising applications to detecting symplectic structures on
4--manifolds. These modules are similar to but different from those
considered by the author, Kent Orr and Peter Teichner and are special
cases of the modules considered subsequently by Shelly Harvey for
arbitrary 3--manifolds.
\end{abstract}

\asciiabstract{%
The classical abelian invariants of a knot are the Alexander module,
which is the first homology group of the the unique infinite cyclic
covering space of S^3-K, considered as a module over the (commutative)
Laurent polynomial ring, and the Blanchfield linking pairing defined
on this module. From the perspective of the knot group, G, these
invariants reflect the structure of G^(1)/G^(2) as a module over
G/G^(1) (here G^(n) is the n-th term of the derived series of
G). Hence any phenomenon associated to G^(2) is invisible to abelian
invariants.  This paper begins the systematic study of invariants
associated to solvable covering spaces of knot exteriors, in
particular the study of what we call the n-th higher-order Alexander
module, G^(n+1)/G^(n+2), considered as a Z[G/G^(n+1)$-module.  We show
that these modules share almost all of the properties of the classical
Alexander module. They are torsion modules with higher-order Alexander
polynomials whose degrees give lower bounds for the knot genus.  The
modules have presentation matrices derived either from a group
presentation or from a Seifert surface. They admit higher-order
linking forms exhibiting self-duality. There are applications to
estimating knot genus and to detecting fibered, prime and alternating
knots. There are also surprising applications to detecting symplectic
structures on 4-manifolds. These modules are similar to but different
from those considered by the author, Kent Orr and Peter Teichner and
are special cases of the modules considered subsequently by Shelly
Harvey for arbitrary 3-manifolds.}

\maketitle

\section{Introduction}\label{intro} The success of algebraic
topology in classical knot theory has been largely confined
to {\it abelian} invariants, that is to say to invariants
associated to the unique regular covering space of $\sk$ with
$\bbz$ as its group of covering translations. These
invariants are the {\it classical Alexander module}, which is
the first homology group of this cover considered as a module
over the {\it commutative} ring $\bbz[t,t^{-1}]$, and the {\it
classical Blanchfield linking pairing}. In turn these
determine the {\it Alexander polynomial} and {\it Alexander
ideals} as well as various numerical invariants associated to
the {\it finite} cyclic covering spaces. From the perspective
of the {\it knot group}, $G=\pi_1(\sk)$, these invariants
reflect the structure of $\gs1/\gs2$ as a module over $G/\gs1$
(here $\gs0=G$ and $\gs n=[\gns-,\gns-]$ is the {\it derived
series of} $G$). Hence any phenomenon associated to
$\gs2$ is {\it invisible} to abelian invariants. This paper
attempts to remedy this deficiency by beginning the
systematic study of invariants associated to {\it solvable}
covering spaces of $\sk$, in particular the study of the {\it
higher-order Alexander module}, $\gs n/\gns+$, considered as a
$\bbz[G/\gs n]$--module. Certainly such modules have been
considered earlier but the difficulties of working with
modules over non-commutative, non-Noetherian, non UFD's seems
to have obstructed progress.

Surprisingly, we show that these {\it higher-order Alexander
modules} share most of the properties of the classical Alexander
module. Despite the difficulties of working with modules over
non-commutative rings, there are applications to estimating knot
genus, detecting fibered, prime and alternating knots as well as
to knot concordance. Most of these properties are not restricted
to the derived series, but apply to other series. For simplicity
this greater generality is discussed only briefly herein.

Similar modules were studied in \cite{COT1} \cite{COT2} \cite{CT}
where important applications to knot concordance were achieved.
The foundational ideas of this paper, as well as the tools
necessary to begin it, were already present in \cite{COT1} and for
that I am greatly indebted to my co-authors Peter Teichner and
Kent Orr. Generalizing our work on knots, Shelly Harvey has
studied similar modules for arbitrary $3$--manifolds and has found
several striking applications: lower bounds for the Thurston norm
of a $2$--dimensional homology class that are much better than C.
McMullen's lower bound using the Alexander norm; and new algebraic
obstructions to a $4$--manifold of the form $M^3\x S^1$ admitting a
symplectic structure \cite{Ha}.

Some notable earlier successes in the area of {\it
non-abelian} knot invariants were the Jones polynomial,
Casson's invariant and the Kontsevitch integral. More in the
spirit of the present approach have been the ``metabelian''
{\it Casson--Gordon invariants} and the {\it twisted Alexander
polynomials} of X.S. Lin and P. Kirk and C. Livingston
\cite{KL}. Most of these detect noncommutativity by studying
representations into known matrix groups over {\it
commutative} rings. The relationship (if any) between our
invariants and these others, is not clear at this time.

Our major results are as follows. For any $n\ge0$ there are
torsion modules $\SA^\bbz_n(K)$ and $\SA_n(K)$, whose
isomorphism types are knot invariants, generalizing the
classical integral and ``rational'' Alexander module ($n=0$)
(Sections~\ref{defs}, \ref{torsion}, \ref{rational}).
$\SA_n(K)$ is a finitely generated module over a
non-commutative principal ideal domain $\bbk_n[t^{\pm1}]$
which is a skew Laurent polynomial ring with coefficients in
a certain skew field (division ring) $\bbk_n$. There are {\it
higher-order Alexander polynomials}
$\Delta_n(t)\in\bbk_n[t^{\pm1}]$ (Section~\ref{polynomials}).
If $K$ does not have (classical) Alexander polynomial 1 then
all of its higher modules are non-trivial and
$\Delta_n\neq1$. The degrees $\d_n$ of these higher order
Alexander polynomials are knot invariants and (using some
work of S. Harvey) we show that they give lower bounds for
knot genera which are provably sharper than the classical
bound ($\d_0\le2\genus(K)$) (see Section~\ref{examples}).

\begin{5.4}\label{5.4}
If $K$ is a non-trivial knot and $n\ge1$ then
$\d_0(K)\le\d_1(K)+1\le\d_2(K)+1\le\dots\le\d_n(K)+1\dots\le2\genus(K)$.
\end{5.4}

\begin{7.4} If $K$ is a knot whose (classical) Alexander
polynomial is not 1 and $k$ is a positive integer then there
exists a hyperbolic knot $K_*$, with the same classical
Alexander module as $K$, for which
$\d_0(K_*)<\d_1(K_*)<\dots<\d_k(K_*)$.
\end{7.4}

There exist presentation matrices for these modules obtained
by pushing loops off of a Seifert matrix
(Section~\ref{seifert}). There also exist presentation
matrices obtained from any presentation of the knot group via
free differential calculus (Section~\ref{fox}).There are
higher order bordism invariants, $\b_n$, generalizing the Arf
invariant (Section~\ref{bordism}) and higher order signature
invariants, $\rho_n$, defined using traces on Von Neumann
algebras (Section~\ref{signatures}). These can be used to
detect chirality. Examples are given wherein these are used to
distinguish knots which cannot be distinguished even by the
$\d_n$. There are also higher order linking forms on
$\SA_n(K)$ whose non-singularity exhibits a self-duality in
the $\SA_n(K)$ (Section~\ref{non-singular}).

The invariants $\SA^\bbz_i$, $\d_i$ and $\rho_i$ have very
special behavior on fibered knots and hence give many new
realizable algebraic obstructions to a knot's being fibered
(Section~\ref{fiber}). Moreover using some deep work of P.
Kronheimer and T. Mrowka \cite{Kr2} the $\d_i$ actually give
new algebraic obstructions to the existence of a symplectic
structure on $4$--manifolds of the form $S^1\x M_K$ where
$M_K$ is the zero-framed surgery on $K$. These obstructions
can be non-trivial even when the Seiberg--Witten invariants
are inconclusive!

\begin{9.5} Suppose $K$ is a non-trivial knot. If $K$ is
fibered then all the inequalities in the above Theorem are
equalities. The same conclusion holds if $S^1\x M_K$ admits a
symplectic structure.
\end{9.5}

Section~\ref{fiber} establishes that, given any $n>0$, there
exist knots with $\d_i+1=\d_0$ for $i<n$ but $\d_n+1\neq\d_0$.

The modules studied herein are closely related to the modules
studied in \cite{COT1} \cite{COT2} \cite{CT}, but are
different. In particular for $n>0$ our $\SA_n$ and $\d_n$
have no known special behavior under concordance of knots.
This is because the $\mathcal{A}_n$ reflect only the
fundamental group of the knot exterior, whereas the modules
of \cite{COT1} reflect the fundamental groups of all possible
slice disk exteriors. To further detail the properties of the
higher-order modules of \cite{COT1} (for example their
presentation in terms of a Seifert surface and their special
nature for slice knots) will require a separate paper
although many of the techniques of this paper will carry over.

\section{Definitions of the higher-order Alexander
modules}\label{defs} The classical Alexander modules of a
knot or link or, more generally, of a $3$--manifold are
associated to the first homology of the universal abelian
cover of the relevant $3$--manifold. We investigate the
homology modules of other regular covering spaces canonically
associated to the knot (or $3$--manifold).

Suppose $M_\G$ is a regular covering space of a connected
CW-complex $M$ such that the group $\G$ is identified with a
subgroup of the group of deck (covering) translations. Then
$H_1(M_\G)$ as a $\bbz\G$--module can be called a {\it
higher-order Alexander module}. In the important special case
that $M_\G$ is connected and $\G$ is the full group of
covering transformations, this can also be phrased easily in
terms of $G=\pi_1(M)$ as follows. If $H$ is any normal
subgroup of $G$ then the action of $G$ on $H$ by conjugation
($h\lra g^{-1}hg$) induces a right $\bbz[G/H]$--module
structure on $H/[H,H]$. If $H$ is a characteristic subgroup of
$G$ then the {\it isomorphism type} (in the sense defined
below) of this module depends only on the isomorphism type
of $G$.

The primary focus of this paper will be the case that $M$ is
a classical knot exterior $\sk$ and on the modules arising
from the family of characteristic subgroups known as the {\it
derived series} of $G$ (defined in Section~\ref{intro}).

\begin{defn}\label{module} The {\it $n^{\rm th}$
(integral) higher-order Alexander module}, $\SA^\bbz_n(K)$,
$n\ge0$, of a knot $K$ is the first (integral) homology group
of the covering space of $\sk$ corresponding to $\gns+$,
considered as a right $\bbz[G/\gns+]$--module, i.e.\
$\gns+/G^{(n+2)}$ as a right module over $\bbz[G/\gns+]$.
\end{defn}

Clearly this coincides with the classical (integral)
Alexander module when $n=0$ and otherwise will be called a
{\it higher-order Alexander module}. It is unlikely that
these modules are finitely generated. However S. Harvey has
observed that they are the torsion submodules of the finitely
presented modules obtained by taking homology relative to the
inverse image of a basepoint \cite{Ha}. The analogues of the
classical {\it rational} Alexander module will be discussed
later in Section~\ref{rational}. These {\it are} finitely
generated.

Note that the modules for different knots (or modules for a
fixed knot with different basepoint for $\pi_1$) are modules
over different (albeit sometimes isomorphic) rings. This
subtlety is even an issue for the classical Alexander module.
If $M$ is an $R$--module and $M'$ is an $R'$--module, we say $M$
{\it is (weakly) isomorphic to} $M'$ if there exists a ring
isomorphism $f\co R\to R'$ such that $M$ is isomorphic to $M'$
as $R$--modules where $M'$ is viewed as an $R$--module via
$f$. If $R$ and $R'$ are group rings (or functorially
associated to groups $G$, $G'$) then we say $M$ is {\it
isomorphic} to $M'$ if there is a group isomorphism $g\co G\lra
G'$ inducing a weak isomorphism.

\begin{prop}\label{invariance} If $K$ and $K'$ are equivalent
knots then $\SA^\bbz_n(K)$ is isomorphic to $\SA^\bbz_n(K')$
for all $n\ge0$.
\end{prop}

\begin{proof}[Proof of \ref{invariance}] If $K$ and $K'$ are
equivalent then their groups are isomorphic. It follows that
their derived modules are isomorphic.
\end{proof}

Thus a knot, its mirror-image and its reverse have isomorphic
modules. In order to take advantage of the peripheral
structure, one needs to use the presence of this extra
structure to restrict the class of allowable ring
isomorphisms. This may be taken up in a later paper. However
in Section~\ref{bordism} and Section~\ref{signatures}
respectively we introduce higher-order bordism and signature
invariants which {\it do} use the orientation of the knot
exterior and hence can distinguish some knots from their
mirror images.

\begin{ex}\label{stabilization} If $K$ is a knot whose
classical Alexander polynomial is 1, then it is well known
that its classical Alexander module $\gs1/\gs2$ is zero. But
if $\gs1=\gs2$ then $\gs n=\gns+$ for all $n\ge1$. Thus each
of the higher-order Alexander modules $\SA^\bbz_n$ is also
trivial. Hence these methods do not seem to give new
information on Alexander polynomial 1 knots. However, it is
shown in Corollary~\ref{non-triviality} that if the classical
Alexander polynomial is {\it not} 1, then {\it all} the
higher-order modules are {\it non-trivial}.
\end{ex}

\begin{ex}\label{trefoil} Suppose $K$ is the right-handed
trefoil, $X=\sk$ and $G=\pi_1(X)$. Since $K$ is a fibered
knot we may assume that $X$ is the mapping torus of the
homeomorphism $f\co\Sigma\to\Sigma$ where $\Sigma$ is a
punctured torus and we may assume $f$ fixes $\p\Sigma$
pointwise. Then $\pi_1(\Sigma)=F\<x,y\>$. Let $X_n$ denote the
covering space of $X$ such that $\pi_1(X_n)\cong\gns+$ and
$\SA^\bbz_n(K)=H_1(X_n)$ as a $\bbz[G/\gns+]$ module. Note
that the infinite cyclic cover $X_0$ is homeomorphic to
$\Sigma\x\bbr$ so that $\pi_1(X_0)\cong\gs1\cong F$. Thus
$X_n$ is a regular covering space of $X_0$ with deck
translations $\gs1/\gns+=F/F^{(n)}$. Since
$\pi_1(X_n)=F^{(n)}$, $H_1(X_n)=F^{(n)}/F^{(n+1)}$ as a module
over $\bbz[F/F^{(n)}]$. Therefore if one considers
$\SA^\bbz_n(K)$ as a module over the subring
$\bbz[\gs1/\gns+]=\bbz[F/F^{(n)}]\sbq\bbz[G/\gns+]$ then it is
merely $F^{(n)}/F^{(n+1)}$ as a module over $\bbz[F/F^{(n)}]$
(a module which depends only on $n$ and the rank of the free
group). More topologically we observe that $X_0$ is homotopy
equivalent to the wedge $W$ of 2 circles and $X_n$ is
(homotopy equivalent to) the result of taking $n$ iterated
universal abelian covers of $W$. Let us consider the case
$n=1$ in more detail. Here $X_1$ is homotopy equivalent to
$W_\infty$, as shown in Figure~\ref{lattice}.

\begin{figure}[ht!]
\cl{\relabelbox\small
\epsfxsize 2in \epsfbox{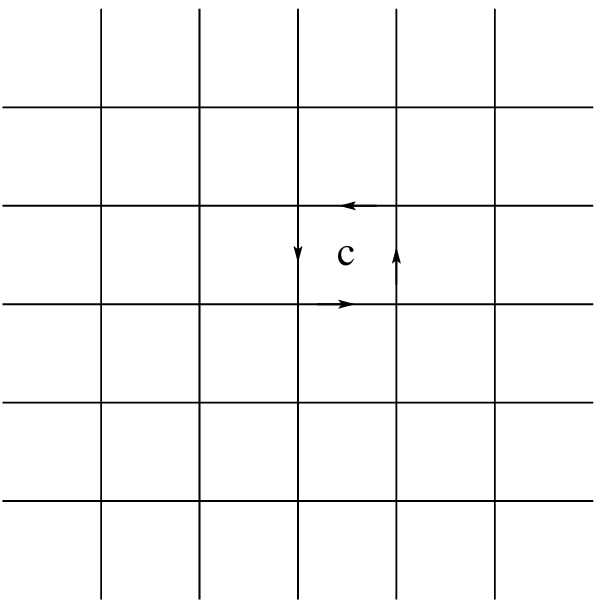}
\adjustrelabel <-2pt, -2pt> {c}{$C$}
\endrelabelbox}
\caption{$W_\infty$} \label{lattice}
\end{figure}

The action of the deck translations
$F/F^{(1)}\cong\bbz\x\bbz$ is the obvious one where $x_*$
acts by horizontal translation and $y_*$ acts by vertical
translation. Clearly $H_1(X_1)$ is an infintely generated
abelian group but as a $\bbz[x^{\pm1},y^{\pm1}]$--module is
cyclic, generated by the loop $C$ in Figure~\ref{lattice}
which represents $xyx^{-1}y^{-1}$ under the identification
$H_1(X_1)\cong F^{(1)}/F^{(2)}$. In fact $H_1(X_1)$ is a free
$\bbz[x^{\pm1},y^{\pm1}]$--module generated by $C$. But
$\SA^\bbz_1(K)=H_1(X_1)$ is a $\bbz[G/\gs2]$--module and so
far all we have discussed is the action of the subring
$\bbz[F/F^{(1)}]=\bbz[\gs1/\gs2]$ because we have completely
ignored the fact that $X_0$ itself has a $\bbz$--action on it.
In fact, since
$1\lra\gs1/\gs2\overset{i}{\lra}G/\gs2\overset{\pi}{\lra}G/\gs1\equiv\bbz\lra1$
is exact, any element of $G/\gs2$ can be written as $gt^m$ for
some $g\in\gs1/\gs2$ and $m\in\bbz$ where $\pi(t)=1$. Thus we
need only specify how $t_*$ acts on $H_1(X_1)$ to describe
our module $\SA^\bbz_1(K)$. To see this action topologically,
recall that, while $X_0$ is homotopy equivalent to $W$, a
more precise description of it is as a countably infinite
number of copies of
$\Sigma\x[-1,1]$ where $\Sigma\x\{1\}\hra(\Sigma\x[-1,1])_i$ is
glued to $\Sigma\x\{-1\}\hra(\Sigma\x[-1,1])_{i+1}$ by the
homeomorphism $f$. Correspondingly, $X_1$ is homotopy
equivalent to $\coprod^\infty_{i=-\infty}(W_\infty\x[-1,1])$
glued together in just such a fashion by lifts of $f$ to
$W_\infty$. Hence $t_*$ acts as $f_*$ acts on
$H_1(X_1)=F^{(1)}/F^{(2)}$. For example if
$f_*(C)=f(xyx^{-1}y^{-1})=w(x,y)C$ then $\SA^\bbz_1(K)$ is a
cyclic module, generated by $C$, with relation
$(t-w(x,y))C=0$. Since $xyx^{-1}y^{-1}$ is represented by the
circle $\p\Sigma$, and since $f$ fixes this circle, in this
case we have that $w(x,y)=1$ and
$\SA^\bbz_1(K)\cong\bbz[G/\gs2]/(t-1)\bbz[G/\gs2]$. This is
interesting because it has $t-1$ torsion represented by the
longitude, whereas the classical Alexander module has no
$t-1$ torsion. This reflects the fact that the longitude
commutes with the meridian as well as the fact that the
longitude, while trivial in $G/\gs2$, is non-trivial in
$\gs2/\gs3\equiv\SA^\bbz_1$.

Since the figure 8 knot is also a fibered genus 1 knot, its
module has a similar form. But note that these modules are not
isomorphic because they are modules over non-isomorphic rings
(since the two knots do not have isomorphic classical
Alexander modules $\gs1/\gs2$). This underscores that the
higher Alexander modules $\SA_i$ should only be used to
distinguish knots with isomorphic $\SA_0,\dots,\SA_{i-1}$.
\qed
\end{ex}

The group of deck translations, $G/\gs n$ of the $\gs n$
cover of a knot complement is solvable but actually satisfies
the following slightly stronger property.

\begin{defn}\label{ptfa} A group $\G$ is poly-(torsion-free
abelian) (henceforth abbreviated PTFA) if it admits a normal
series $\<1\>=G_n\ltri G_{n-1}\ltri\dots\ltri G_0=\G$ such that
the factors $G_i/G_{i+1}$ are torsion-free abelian (Warning -
in the group theory literature only a subnormal series is
required).
\end{defn}

This is a convenient class (as we shall see) because it is
contained in the class of {\it locally indicable} groups
\cite[Proposition 1.9]{Str} and hence $\bbz\G$ is an integral
domain \cite{Hig}. Moreover it is contained in the class of {\it
amenable} groups and thus $\bbz\G$ embeds in a classical quotient
(skew) field \cite[Theorem 5.4]{Do}.

It is easy to see that every PTFA group is solvable and
torsion-free and although the converse is not quite true,
every solvable group such that each $\gs n/\gns+$ is
torsion-free, is PTFA. Every torsion-free nilpotent group is
PTFA.

Consider a tower of regular covering spaces
$$
M_n\lra M_{n-1}\lra\dots\lra M_1\lra M_0 = M
$$
such that each $M_{i+1}\lra M_i$ has a torsion-free abelian
group of deck translations and each $M_i\lra M$ is a regular
cover. Then the group $\G$ of deck translations of $M_n\lra
M$ is PTFA and it is easy to see that such towers correspond
precisely to normal series for such a group.

\begin{ex} If $G=\pi_1(\sk)$ and $\gs n$ is the
$n^{\rm th}$ term of the derived series then $G/\gs n$ is
PTFA since each $\gs i/G^{(i+1)}$ is known to be torsion free
\cite{Str}. Therefore taking iterated universal abelian
covers of $S^3-K$ yields a PTFA tower as above. Hence the
$n^{\rm th}$ higher-order Alexander module generalizes the
classical Alexander module in that the latter is the case of
taking a single universal abelian covering space.
\end{ex}

There is certainly more information to be found in modules
obtained from {\it other} $\G$--covers. For most of the proofs
we can consider a general $\G$--cover where $\G$ is PTFA. Thus
there are other families of subgroups which merit scrutiny,
and are covered by most of the theorems to follow, but which
will not be discussed in this paper. Primary among these is
the lower central series of the commutator subgroup of $G$.

For a general $3$--manifold with first Betti number equal to
$1$ (which we cover since it is no more difficult than a knot
exterior) it is necessary to use the {\it rational derived
series} to avoid zero divisors in the group ring:

\begin{ex} For {\it any} group $G$, the $n^{\rm th}$ term
of the {\it rational derived series} is defined by
$\gs0_\BQ=G$ and
$\gs n_\BQ=[G^{(n-1)}_\BQ,G^{(n-1)}_\BQ]\cd N$ where
$N=\{g\in G^{(n-1)}_\BQ|$ some non-zero power of $g$ lies in
$[G^{n-1}_\BQ,G^{n-1}_\BQ]\}$. It is easy to see that
$G/\gs n_\BQ$ is PTFA. This corresponds to taking iterated
universal {\it torsion-free} abelian covering spaces. For
knot groups, $\gs n_\BQ=\gs n$ \cite{Str}.
\end{ex}

\begin{defn}\label{general def} If $M$ is an arbitrary
connected CW-complex with fundamental group $G$, then the {\it
$n^{\rm th}$ (integral) higher-order Alexander} module,
$\SA^\bbz_n(M)$,\linebreak
$n\ge0$, of $M$ is $H_1(M_n;\bbz)$ ($M_n$ is
the cover of $M$ with $\pi_1(M_n)=\gns+_\bbq)$ considered as a
right $\bbz[G/\gns+_\bbq]$--module.
\end{defn}

\sh{More on the relationship of $\SA^\bbz_n(K)$ to
$\pi_1(\sk)$}

We have seen that if $H$ is any characteristic subgroup of
$G$ then the isomorphism type of $H/[H,H]$, as a right module
over $\bbz[G/H]$, is an invariant of the isomorphism type of
$G$. Moreover, $\SA^\bbz_n(K)$ has been defined as this
module in the case $G=\pi_1(\sk)$ and $H=\gns+$. The
following elementary observation clarifies this relationship.
Its proof is left to the reader. One consequence will be that
for any knot there exists a hyperbolic knot with isomorphic
$\SA^\bbz_n$ for all $n$.

\begin{prop}\label{derived} Suppose $f\co G\lra P$ is an
epimorphism. Then $f$ induces isomorphisms
$f_n\co\SA^\bbz_n(G)\lra\SA^\bbz_n(P)$ for all $n\le m$ if and
only if the kernel of $f$ is contained in $G^{(m+2)}_\bbq$.
Hence $f$ induces such isomorphisms for all finite $n$ if and
only if kernel $f\subset\bigcap^\infty_{n=1}\gs n_\bbq$.
\end{prop}

\begin{cor}\label{hyperbolic} For any knot $K$, there is a
hyperbolic knot $\wt K$ and a degree one map
$f\co S^3\backslash\wt K\lra\sk$ (rel boundary) which induces
isomorphisms $\SA^\bbz_n(\wt K)\lra\SA^\bbz_n(K)$ for all $n$.
\end{cor}

\begin{proof}[Proof of Corollary~\ref{hyperbolic}] In fact it
is known that $\wt K$ can be chosen so that the kernel of
$f_*$ is a perfect group (or in other words that $f$ induces
isomorphisms on homology with $\bbz[\pi_1(\sk)]$
coefficients). The first reference I know to this fact is by
use of the ``almost identical link imitations'' of Akio
Kawauchi \cite[Theorem 2.1 and Corollary 2.2]{Ka}. A more
recent and elementary construction can be adopted from
\cite[Section 4]{BW}. Any perfect subgroup is contained in its
own commutator subgroup and hence, by induction, lies in
every term of the derived series. An application of
Proposition~\ref{derived} finishes the proof.

\begin{ex}
If $K'$ is a knot and $K$ is a knot
whose (classical) Alexander polynomial is $1$ then $K'$ and
$K'\#K$ have isomorphic higher-order modules since there is a
degree one map $S^3\backslash(K'\#K)\to S^3\backslash K'$
which induces an epimorphism on $\pi_1$ whose kernel is
$\pi_1(\sk)^{(1)}$. The observation then follows from
Proposition~\ref{derived} and Example~\ref{stabilization}.
\end{ex}

\end{proof}

\section{Properties of higher-order Alexander modules of
knots: Torsion}\label{torsion}

In this section we will show that higher-order Alexander
modules have one key property in common with the classical
Alexander module, namely they are torsion-modules. In
Section~\ref{non-singular} we define a linking pairing on
these modules which generalizes the Blanchfield linking
pairing on the Alexander module. All of the results of this
section follow immediately from \cite[Section 2]{COT1} but a
simpler proof of the main theorem is given here.

A right module $A$ over a ring $R$ is said to be a {\it
torsion module\/} if, for any $a\in A$, there exists a
non-zero-divisor $r\in R$ such that $ar=0$.

Our first goal is:

\begin{thm}\label{Torsion} The higher-order Alexander modules
$\SA^\bbz_n(K)$ of a knot are torsion modules.
\end{thm}

This is a consequence of the more general result which
applies to any complex $X$ with $\pi_1(X)$ finitely-generated
and $\b_1(X)=1$ and any PTFA $\G$
\cite[Proposition~2.11]{COT1} but we shall give a different,
self-contained proof (Proposition~\ref{rank}). The more
general result will be used in later chapters to study
general $3$--manifolds with $\b_1=1$.

Suppose $\G$ is a PTFA group. Then $\bbz\G$ has several
convenient properties --- it is an integral domain and it has
a classical field of fractions. Details follow.

Recall that if $A$ is a {\em commutative\/} ring and $S$ is a
subset closed under multiplication, one can construct the
{\em ring of fractions\/} $AS^{-1}$ of elements $as^{-1}$
which add and multiply as normal fractions. If $S=A-\{0\}$ and
$A$ has no zero divisors, then $AS^{-1}$ is called the {\em
quotient field\/} of $A$. However, if $A$ is {\em
non-commutative\/} then $AS^{-1}$ does not always exist (and
$AS^{-1}$ is not a priori isomorphic to $S^{-1}\ngth A$). It
is known that if $S$ is a {\em right divisor set\/} then
$AS^{-1}$ exists (~\cite[p.~146]{P} or ~\cite[p.~52]{Ste}).
If $A$ has no zero divisors and $S=A-\{0\}$ is a right divisor
set  then $A$ is called an {\em Ore domain\/}. In this case
$AS^{-1}$ is a skew field, called the {\em classical right
ring of quotients\/} of $A$. We will often refer to this
merely as the {\em quotient field } of $A$ . A good reference
for non-commutative rings of fractions is Chapter~2 of
\cite{Ste}. In this paper we will always use {\em right\/}
rings of fractions.

\begin{prop}\label{ore} If $\G$ is PTFA then $\bbq\G$ (and
hence $\bbz\G$) is a right (and left) Ore domain; i.e.\
$\bbq\G$ embeds in its classical right ring of quotients
$\SK$, which is a skew field.
\end{prop}

\begin{proof} For the fact (due to A.A. Bovdi) that $\bbz\G$
has no zero divisors see \cite[pp.~591--592]{P} or
\cite[p.~315]{Str}. As we have remarked, any PTFA group is
solvable. It is a result of J.~Lewin~\cite{Lew} that for
solvable groups such that $\bbq\G$ has no zero divisors,
$\bbq\G$ is an Ore domain (see Lemma~3.6 iii p.~611 of
\cite{P}). It follows that $\bbz\G$ is also an Ore domain.
\end{proof}

\begin{rem}\label{facts} Skew fields share many of the key
features of (commutative) fields. We shall need the following
elementary facts about the right skew field of quotients
$\SK$. It is naturally a $\SK$--$\SK$--bimodule and a
$\bbz\G$--$\bbz\G$--bimodule.
\begin{description}
\item[Fact 1]  $\SK$ is flat as a left $\bbz\G$--module, i.e.\
$\cdot\ox_{\bbz\G}\SK$ is
exact~\cite[Proposition~II.3.5]{Ste}.
\item[Fact 2] Every module over $\SK$ is a free
module~\cite[Proposition~I.2.3]{Ste} and such modules have a
well defined rank $\rks$ which is additive on short exact
sequences ~\cite[~p.~48]{Co2}.
\end{description}
\end{rem}

If $\SA$ is a module over the Ore domain $\SR$ then the {\it
rank of $\SA$} denotes $\rank_\SK(\SA\ox_\SR\SK)$. $\SA$ is a
torsion module if and only if $\SA\ox_{\SR}\SK=0$ where $\SK$
is the quotient field of $\SR$, i.e.\ if and only if the rank
of $\SA$ is zero \cite[II Corollary~3.3]{Ste}. In general,
the set of torsion elements of $A$ is a submodule which is
characterized as the kernel of $\SA\to\SA\ox_{\SR}\SK$. Note
that if $\SA\cong\SR^r\op$(torsion) then $\rank\SA=r$.

\begin{description}
\item[Fact 3] If $\SC$ is a non-negative finite chain complex of
finitely generated free (right) $\bbz\G$--modules then the
equivariant Euler characteristic, $\chi(\SC)$,
 given by $\sum^\infty_{i=0}(-1)^i\rank C_i$, is defined and equal to
$\sum^\infty_{i=0}(-1)^i\rank H_i(\SC)$ and
$\sum^\infty_{i=0}(-1)^i\rank H_i(\SC\ox_{\bbz\G}\SK)$. This is an
elementary consequence of Facts~1 and 2.
\end{description}

There is another especially important property of PTFA
groups (more generally of locally indicable groups) which
should be viewed as a natural generalization of properties of
the free abelian group. This is an algebraic generalization
of the (non-obvious) fact that any infinite cyclic cover of a
$2$--complex with vanishing $H_2$ also has vanishing $H_2$
(see Proposition~\ref{injectivity}).

\begin{prop}\label{strebel} {\rm(R.
Strebel~\cite[p.~305]{Str})} Suppose $\G$ is a PTFA group and
$R$ is a commutative ring. Any map between projective right
$R\G$--modules whose image under the functor $-\ox_{R\G}R$ is
injective, is itself injective.
\end{prop}

We can now offer a simple proof of Theorem~\ref{Torsion}.

\begin{proof}[Proof of Theorem~\ref{Torsion}] The knot
exterior has the homotopy type of a finite connected
$2$--complex $Y$ whose Euler characteristic is $0$. Let
$\G=G/\gns+$ and let $\SC=(0\lra C_2\xrightarrow{\partial_2}
C_1\xrightarrow{\partial_1} C_0\lra 0)$ be the free $\bbz\G$
cellular chain complex for $Y_\G$ (the $\G$--cover of $Y$ such
that $\pi_1(Y)=\gns+$) obtained by lifting the cell structure
of $Y$. Then $\chi(\SC)=\chi(Y)=0$. It follows from Fact~3
that $\rank H_2(Y_\G)-\rank H_1(Y_\G)+\rank H_0(Y_\G)=0$. Now
note that $(\SC,\p)$ is sent, under the augmentation
$\e\co\bbz\G\lra\bbz$, to
($\SC\ox_{\bbz\G}\bbz,\p\ox_{\bbz\G}\id$) which can be
identified with the chain complex for the original cell
structure on $Y$. Since $H_2(Y;\bbz)=0$, $\p_2\ox\id$ is
injective. By Proposition~\ref{strebel}, it follows that
$\p_2$ itself is injective, and hence that $H_2(Y_\G)=0$.

Now we claim that $H_0(Y_\G)$ is a torsion module. This is
easy since $H_0(Y_\G)\cong\bbz$. If $H_0(Y_\G)$ were {\it not}
torsion then $1\in\bbz$ generates a free $\bbz\G$ submodule.
Note that $\G$ is not trivial since $G\neq\gs1$. This is a
contradiction since, as an abelian group, $\bbz\G$ is free on
more than one generator and hence cannot be a subgroup of
$\bbz$.
\end{proof}

Now that we have proved that the higher-order modules of a
knot are torsion modules, we look at the homology of covering
spaces in more detail and in a more abstract way. This point
of view allows for greater generality and for more concise
notation. Viewing homology of covering spaces as homology
with twisted coefficients clarifies the calculations of the
homology of induced covers over subspaces.

\medskip
\sh{Homology of PTFA covering spaces}

Suppose $X$ has the homotopy type of a connected CW-complex,
$\G$ is any group and $\phi\co\pi_1(X,x_0)\lra\G$ is a
homomorphism. Let $X_\G$ denote the {\em regular $\G$--cover
of $X$ associated to $\phi$} (by pulling back the universal
cover of $B\G$ viewed as a principal $\G$--bundle). If $\phi$
is surjective then $X_\G$ is merely the connected covering
space $X$ associated to $\Ker(\phi)$. Then $X_\G$ becomes a
{\it right} $\G$--set as follows. Choose a point $*\in
p^{-1}(x_0)$. Given $\g\in\G$, choose a loop $w$ in $X$ such
that $\phi([w])=\g$. Let $\tl w$ be a lift of $w$ to
$X_\G$ such that $\tl w(0)=*$. Let $d_w$ be the unique
covering translation such that $d_w(*)=\tl w(1)$. Then $\g$
acts on $X_\G$ by $d_w$. This merely the ``usual'' left
action \cite[Section 81]{M2}. However, for certain historical
reasons we shall use the associated {\em right} action where
$\g$ acts by $(d_w)^{-1}$. If $\phi$ is not surjective and we
set $\pi=\mbox{image}(\phi)$ then $X_\G$ is a disjoint union
of copies of the connected cover $X_\pi$ associated to
$\Ker(\phi)$. The set of copies is in bijection with the set
of right cosets $\G/\pi$. In fact it is best to think of
$p^{-1}(x_0)$ as being identified with $\G$. Then $\G$ acts
on $p^{-1}(x_0)$ by right multiplication. If $\g\in\pi$, then
$\g$ sends $*$ to the endpoint of the path $\tl w$ such that
$\tl w(0)=*$ and $\phi([w])=\g^{-1}$. Hence $*$ and $(*)\g$
are in the same path component of $X_\G$. If $\tau\in\G$ is a
non-trivial coset representative then $(*)\tau$ lies in a
different path component than $*$. But the path $\tl w$,
acted on by the deck translation corresponding to $\tau$,
begins at $(*)\tau$ and ends at $(\tl
w(1))\tau=(*)(\g)(\tau)=(*)(\g\tau)$. Thus
$(*)\tau$ and $(*)\tau'$ lie in the same path component if
and only if they lie in the same right coset $\pi\tau$ of
$\G/\pi$.

For simplicity, the following are stated for the ring $\bbz$,
but also hold for $\bbq$. Let $\SM$ be a $\bbz\G$--bimodule
(for us usually $\bbz\G$, $\SK$, or a ring $\SR$ such that
$\bbz\G\subset\SR\subset\SK$, or $\SK/\SR$). The following
are often called the equivariant homology and cohomology of
$X$.

\begin{defn}\label{homology} Given $X$, $\phi$, $\SM$ as
above, let
$$
H_{\ast}(X;\SM)\equiv H_{\ast}(\SC(X_\G;\bbz)\ox_{\bbz\G}\SM)
$$
as a right $\bbz\G$ module, and $H^{\ast}(X;\SM)\equiv
H_{\ast}\left(\Hom_{\bbz\G}(\SC(X_\G;\bbz),\SM)\right)$ as a
left $\bbz\G$--module.
\end{defn}

These are also well-known to be isomorphic (respectively) to
the homology (and cohomology) of $X$ with coefficient system
induced by $\phi$ (see Theorems~VI 3.4 and 3.4$^*$ of
\cite{W}). The advantage of this formulation is that it
becomes clear that the surjectivity of $\phi$ is irrelevant.

\begin{rem}\ 

\begin{enumerate}
\item  Note that $H_{\ast}(X;\bbz\G)$ as in
Definition~\ref{homology} is merely
$H_{\ast}(X_\G;\bbz)$ as a right $\bbz\G$--module. Thus
$\SA^\bbz_n\cong H_1(\sk;\bbz\G)$ where $\G=G/\gns+$ and
$G=\pi_1(\sk)$. Moreover if $\SM$ is flat as a left
$\bbz\G$--module then
$H_{\ast}(X;\SM)\cong H_{\ast}(X_\G;\bbz)\ox_{\bbz\G}\SM$. In
particular this holds for $\SM=\SK$ by 3.3. Thus
$H_{\ast}(X_\G)=H_{\ast}(X;\bbz\G)$ is a torsion module if and
only if $H_{\ast}(X;\SK)=H_{\ast}(X_\G)\ox_{\bbz\G}\SK=0$ by
the remarks below 3.3.
\item  Recall that if $X$ is a compact, oriented $n$--manifold
then by Poincar\'e duality $H_p(X;\SM)$ is isomorphic to
$H^{n-p}(X,\partial X;\SM)$ which is made into a right
$\bbz\G$--module using the obvious involution on this group
ring \cite{Wa}.
\item We also have a universal coefficient spectral sequence
as in \cite[Theorem 2.3]{L3}. This collapses to the usual
Universal Coefficient Theorem for coefficients in a
(noncommutative) principal ideal domain (in particular for
the skew field $\SK$). Hence
$H^n(X;\SK)\cong\Hom_\SK(H_n(X;\SK),\SK)$. In this paper we
only need the UCSS in these special cases where it coincides
with the usual UCT.
\end{enumerate}
\end{rem}

We now restrict to the case that $\G$ is a PTFA group and
$\SK$ is its (skew) field of quotients. We investigate $H_0$,
$H_1$ and $H_2$ of spaces with coefficients in $\bbz\G$ or
$\SK$.

\begin{prop}\label{Ho} Suppose $X$ is a connected CW complex.
If $\phi\co\pi_1(X)\lra\G$ is a non-trivial coefficient system
then $H_0(X;\SK)=0$ and $H_0(X;\bbz\G)$ is a torsion module.
\end{prop}

\begin{proof} By \cite[p. 275]{W} and \cite[p.34]{Br},
$H_0(X;\SK)$ is isomorphic to the cofixed set $\SK/\SK I$
where $I$ is the augmentation ideal of $\bbz\pi_1(X)$ acting
via $\pi_1(X)\lra\G\lra\SK$. If $\phi$ is non-zero then this
composition is non-zero and hence $I$ contains an element
which acts as a unit. Hence $\SK I=\SK$.
\end{proof}

The following lemma summarizes the basic topological
application of Strebel's result (Proposition~\ref{strebel}).

\begin{prop}\label{injectivity} Suppose $(Y,A)$ is a
connected $2$--complex with $H_2(Y,A;\bbq)\linebreak
\cong0$ and suppose
$\phi\co\pi_1(Y)\lra\G$ defines a coefficient system on $Y$ and
$A$ where $\G$ is a PTFA group. Then $H_2(Y,A;\bbz\G)=0$, and
so $H_1(A;\bbz\G)\lra H_1(Y;\bbz\G)$ is injective.
\end{prop}

\begin{proof} Let $\SC$ be the free $\bbz\G$ chain complex
for the cellular structure on $(Y_\G,A_\G)$ (the $\G$--cover
of $Y$) obtained by lifting the cell structure of $(Y,A)$. It
suffices to show $\p_2\co C_2\lra C_1$ is a monomorphism. By
Proposition~\ref{strebel} this will follow from the
injectivity of $\p_2\ox\id\co C_2\ox_{\bbz\G}\bbz\lra
C_1\ox_{\bbz\G}\bbz$. But this map can be canonically
identified with the corresponding boundary map in the
cellular chain complex of $(Y,A)$, which is injective since
$H_2(Y,A;\bbq)\cong H_2(Y,A;\bbz)\cong0$.
\end{proof}

The following lemma generalizes the key argument of the proof
of Theorem~\ref{Torsion}.

\begin{lem}\label{2complex} Suppose $Y$ is a connected
$2$--complex with $H_2(Y;\bbz)=0$ and $\phi\co\pi_1(Y)\lra\G$ is
non-trivial. Then $H_2(Y;\SK)=0$; and if $Y$ is a finite
complex then $\rks H_1(Y;\SK)=\b_1(Y)-1$.
\end{lem}

\begin{proof} By Proposition~\ref{injectivity}
$H_2(Y;\bbz\G)=0$ and $H_2(Y;\SK)=0$ by Remark~3.6.1. Since
$\phi$ is non-trivial, Proposition~\ref{Ho} implies that
$H_0(Y;\SK)=0$. But by Fact 3 (as in the proof of
Theorem~\ref{Torsion}) $\rank_\SK H_2(Y;\SK)-\rank_\SK
H_1(Y;\SK)+\rank_\SK H_0(Y;\SK)=1-\b_1(Y)$ and the result
follows.
\end{proof}

Note that if $\b_1(Y)=0$ then any homomorphism from
$\pi_1(Y)$ to a PTFA group is necessarily the zero
homomorphism.

\begin{prop}\label{rank} Suppose $\pi_1(X)$ is
finitely-generated and $\phi\co\pi_1(X)\lra\G$ is non-trivial.
Then
$$
\rank_\SK H_1(X;\bbz\G)\le\b_1(X) - 1.
$$
In particular, if $\b_1(X)=1$ then $H_1(X;\bbz\G)$ is a
torsion module.
\end{prop}

\begin{proof} Since the first homology of a covering space of
$X$ is functorially determined by $\pi_1(X)=G$, we can
replace $X$ by a $K(G,1)$. We will now construct an
epimorphism $f\co E\lra G$ from a group $E$ which has a very
efficient presentation. Suppose
$H_1(G)\cong\bbz^m\x\bbz_{n_1}\x\dots\x\bbz_{n_k}$. Then
there is a finite generating set
$\{g_1,\dots,g_m,g_{m+1},\dots,g_{m+k},\dots|i\in I\}$ for
$G$ such that $\{g_1,\dots,g_{m+k}\}$ is a ``basis'' for
$H_1(G)$ wherein if $i>m+k$ then $g_i\in[G,G]$ and if $m<i\le
m+k$ then $g^{n_i}_i\in[G,G]$. Consider variables $\{x_j|j\in
I\}$. Hence for each $i$ there is a word $w_i(x_1,\dots)$ in
these variables such that $w_i$ lies in the commutator
subgroup of the free group on $\{x_j\}$, and such that if
$i>m+k$ then $g_i=w_i(g_1,\dots)$ and if
$m<i\le m+k$ then $g^{n_i}_i=w_i(g_1,\dots)$. Let $E$ have
generators $\{x_i|i\in I\}$ and relations $\{x_i=w_i|i>m+k\}$
and $\{x^{n_i}_i=w_i|m<i\le m+k\}$. The obvious epimorphism
$f\co E\lra G$ given by $f(x_i)=g_i$ is an $H_1$--isomorphism.
The composition $\phi\circ f$ defines a $\G$ covering space
of $K(E,1)$. Since $f$ is surjective we can build $K(G,1)$
from $K(E,1)$ by adjoining cells of dimensions at least 2.
Thus $H_1(G,E;\bbz\G)=0$ because there are no relative
$1$--cells and consequently
$f_*\co H_1(E;\bbz\G)\lra H_1(G;\bbz\G)$ is also surjective.
Since $\SK$ is a flat
$\bbz\G$ module $f_*\co H_1(E;\SK)\lra H_1(G;\SK)$ is
surjective. Thus $\rank_\SK H_1(X;\bbz\G)=\rank_\SK
H_1(X;\SK)\le\rank_\SK H_1(E;\SK)$. Now note that
$E=\pi_1(Y)$ where $Y$ is a connected, finite $2$--complex
(associated to the presentation) which has vanishing second
homology. Again since $H_1$ is functorially determined by
$\pi_1$, $H_1(E;\SK)\cong H_1(Y;\SK)$. Lemma~\ref{2complex}
above shows that
$\rank_\SK H_1(Y;\SK)=\b_1(Y)-1=\b_1(E)-1=\b_1(X)-1$ and the
result follows.
\end{proof}

\begin{ex}
It is somewhat remarkable (and turns out to be crucially important) that the
previous two results fail without the finiteness assumption. If
Proposition~\ref{rank} were true without the finiteness assumption, all of
the inequalities of Theorem~\ref{non-decreasing} would be equalities.
Consider $E=\langle x,z_i\mid z_i=[z_{i+1},x],i\in\bbz\rangle$. This is the
fundamental group of an (infinite) $2$--complex with $H_2=0$. Note that
$\b_1(E)=1$. But the abelianization of $E^{(1)}$ has a presentation
$\langle z_i\mid z_i=(1-x)z_{i+1}\rangle$ as a module over $\bbz[x^{\pm1}]$
and thus has rank~1, {\it not} $\b_1(E)-1$ as would be predicted by
Proposition~\ref{rank}.
\end{ex}

\begin{cor}\label{acyclic} Suppose $M$ is a compact, orientable, connected
$3$--manifold such that $\b_1(M)=1$. Suppose
$\phi\co\pi_1(M)\lra\G$ is a homomorphism that is non-trivial on
abelianizations where $\G$ is PTFA. Then $H_*(M,\p M;\SK)\cong0
\cong H_*(M;\SK)$.
\end{cor}

\begin{proof} Propositions \ref{Ho} and \ref{rank} imply
$H_0(M;\SK)\cong H_1(M;\SK)\cong0$. Since it is well known
that the image of $H_1(\p M;\bbq)\lra H_1(M;\bbq)$ has
one-half the rank of $H_1(\p M;\bbq)$, $\p M$ must be either
empty or a torus. Suppose the latter. Then this
inclusion-induced map is surjective. Therefore the induced
coefficient system $\phi\circ i_*\co\pi_1(\p M)\lra\G$ is
non-trivial since it is non-trivial on abelianizations. Thus
$H_0(\p M;\SK)=0$ by Proposition~\ref{Ho}, implying that
$H_1(M,\p M;\SK)=0$. By Remark 3.6, $H_2(M;\SK)\cong
H^1(M,\p M;\SK)\cong\Hom(H_1(M,\p M;\SK),\SK)\cong0$.
Similarly $H_3(M;\SK)\cong0$. Then
$H_*(M;\SK)\cong0\Rightarrow H_*(M,\p M;\SK)\cong0$ by
duality and the universal coefficient theorem.
\end{proof}

Thus we have shown that the definition of the classical
Alexander module, i.e.\ the torsion module associated to the
first homology of the infinite cyclic cover of the knot
complement, can be extended to {\em higher-order Alexander
modules} $\SA^\bbz_\G=H_1(M;\bbz\G)$ which are $\bbz\G$
torsion modules associated to {\em arbitrary\/} PTFA covering
spaces. Indeed, by Proposition~\ref{rank}, this is true for
any finite complex with $\b_1(M)=1$.

\section{Localized higher-order modules}\label{rational}

In studying the classical abelian invariants of knots, one
usual studies not only the ``integral'' Alexander module,
$H_1(\sk;\bbz[t,t^{-1}])$, but also the {\it rational
Alexander module} $H_1(\sk;\bbq[t,t^{-1}])$. Even though some
information is lost in this localization, $\bbq[t,t^{-1}]$ is
a principal ideal domain and one has a good classification
theorem for finitely generated modules over a PID. Moreover
the rational Alexander module is {\it self-dual} whereas the
integral module is not \cite{Go}. In considering the
higher-order modules it is even more important to localize our
rings $\bbz[G/\gs n]$ in order to define a higher-order
``rational'' Alexander module over a (non-commutative) PID.
Here, significant information will be lost but this
simplification is crucial to the definition of numerical
invariants. Recall that an integral domain is a {\it right
(respectively left) PID} if every right (respectively left)
ideal is principal. A ring is a PID if it is both a left and
right PID. The definition of the relevant PID's follows.

Let $G$ be a group with $\b_1(G)=1$ and let $\G_n=G/\gns+_\bbq$
(which is the same as the ordinary derived series for a knot
group). Recall that the (integral) Alexander module was defined as
$\SA^\bbz_n(G)=H_1(G;\bbz\G_n)$ in Definition~\ref{module} and
Definition~\ref{general def}. Below we will describe a PID $R_n$
such that $\bbq\G_n\subset R_n\subset\SK_n$ and such that $R_n$ is
a localization of $\bbq\G_n$, i.e.\ $R_n=\bbq\G_n(S^{-1})$ where
$S$ is a right divisor set in $\bbq\G_n$. Using this we define the
``localized'' derived modules. These will be analyzed further in
Section~\ref{polynomials}. These PID's were crucial in our
previous work \cite{COT1}.

\begin{defn}\label{an} The {\it $n^{\rm th}$ ``localized'' Alexander
module} of a knot $K$, or, simply, the {\it $n^{\rm th}$ Alexander
module} of $K$ is $\SA_n(K)=H_1(\sk;R_n)$.
\end{defn}

\begin{prop}\label{finite generation} The $n^{\rm th}$
Alexander module is a finitely-generated torsion module over
the PID $R_n$.
\end{prop}

\begin{proof} Let $M_n$ denote the covering space of $M=\sk$
with $\pi_1(M_n)=\gns+$. Then $\SA_n(K)$ is the first
homology of the chain complex $C_*(M_n)\ox_{\bbz\G_n}R_n$.
This is a chain complex of finitely generated free
$R_n$--modules since $M$ has the homotopy type of a finite
complex and we can use the lift of this cell structure to
$M_n$. Since a submodule of a finitely-generated free module
over a PID is again a finitely-generated free module
(\cite{J}, Theorem~17), it follows that the homology groups
are finitely generated.
\end{proof}

Now we define the rings $R_n$ and show that they are PID's by
proving that they are isomorphic to {\it skew Laurent
polynomial rings} $\bbk_n[t^{\pm1}]$ over a skew field
$\bbk_n$. This makes the analogy to the classical rational
Alexander module even stronger.

Before defining $R_n$ in general, we do so in a simple
example.

\begin{ex}\label{trefoil2} We continue with
Example~\ref{trefoil} where $G=\pi_1(\sk)$ and $K$ is a
trefoil knot. We illustrate the structure of
$\bbz[G/\gs2]=\bbz\G_1$ as a skew Laurent polynomial ring in
one variable with coefficients in $\bbz[\gs1/\gs2]$. Recall
that since the trefoil knot is fibered, $\gs1/\gs2\cong
F/F^{(1)}\cong\bbz\x\bbz$ generated by $\{x,y\}$. Hence
$\bbz[\gs1/\gs2]$ is merely the (commutative) Laurent
polynomial ring $\bbz[x^{\pm1},y^{\pm1}]$. If we choose, say,
a meridian $\mu\in G/\gs2$ then $G/\gs2$ is a semi-direct
product $\gs1/\gs2\rtimes\bbz$ and any element of $G/\gs2$
has a unique representative $\mu^mg$ for some $m\in\bbz$ and
$g\in\gs1/\gs2$, i.e.\ $\mu^mx^py^q$ for some integers $m$,
$p$, $q$. Thus any element of $\bbz[G/\gs2]$ has a canonical
representation of the form
$\sum^\infty_{m=-\infty}\mu^mp_m(x,y)$ where
$p_m(x,y)\in\bbz[x^{\pm1},y^{\pm1}]$. Hence $\bbz[G/\gs2]$
can be identified with the Laurent polynomial ring in one
variable $\mu$ (or $t$ for historical significance) with
coefficients in the Laurent polynomial ring
$\bbz[x^{\pm1},y^{\pm1}]$. Observe that the product of 2
elements in canonical form is not in canonical form. However,
for example,
$(x^py^q)\cd\mu=\mu(\mu^{-1}x^py^q\mu)=\mu((x^py^q)\mu_*)$.
Hence this is not a true polynomial ring, rather the
multiplication is twisted by the automorphism $\mu_*$ of
$\bbz[\gs1/\gs2]$ induced by conjugation $g\to\mu^{-1}g\mu$
(the action of the generator $t\in\bbz$ in the semi-direct
product structure). The action $\mu_*$ (or $t_*$) is merely
the action of $t$ on the Alexander module of the trefoil
$\bbz[t,t^{-1}]/t^2-t+1\cong\bbz\x\bbz$ with basis $\{x,y\}$.

Moreover this skew polynomial ring
$\bbz[\gs1/\gs2][t^{\pm1}]$ embeds in the ring
$R_1=\bbk_1[t^{\pm1}]$, where $\bbk_1$ is the quotient field
of the coefficient ring $\bbz[x^{\pm1},y^{\pm1}]$ (in this
case the (commutative) field of rational functions in the 2
commuting variables $x$ and $y$). Thus $\bbz[G/\gs2]$ embeds
in this (noncommutative) PID $R_1$ (this is proved below)
that also has the structure of a skew Laurent polynomial ring
over a field. Note that, under this embedding, the subring
$\bbz[\gs1/\gs2]$ is sent into the subring of polynomials of
degree $0$, i.e.\ $\bbk_1$ and this embedding is just the
canonical embedding of a commutative ring into its quotient
field (and is thus independent of the choice of $\mu$!). \qed
\end{ex}

Now we define $R_n$ in general. Let $\wt G_n$, $n\ge1$, be the
kernel of the map $\pi\co G/\gs n_\bbq\lra G/\gs1_\bbq$ (the
latter is infinite cyclic by the hypothesis that $\b_1(G)=1$.
For the important case that $G$ is a knot group, $\wt G_n$ is
the commutator subgroup modulo the $n^{\rm th}$ derived
subgroup. Since $G/\gs n_\bbq$ is PTFA by Example~2.7, the
subgroup $\wt G_n$ is also PTFA. Thus $\bbz[\wt G_n]$ is an
Ore domain by Proposition~\ref{ore}. Let $S_n=\bbz[\wt
G_{n+1}]-\{0\}$, $n\ge0$, a subset of
$\bbz\G_n=\bbz[G/\gns+_\bbq]$. By \cite[p.~609]{P} $S_n$ is a
right divisor set of $\bbz\G_n$ and we set
$R_n=(\bbz\G_n)(S_n)^{-1}$. Hence $\bbz\G_n\sbq
R_n\sbq\SK_n$. Note that $S_0=\bbz-\{0\}$ so $R_0=\bbq[J]$
where $J$ is the infinite cyclic group $G/\gs1_\bbq$,
agreeing with the classical case. By Proposition II.3.5
\cite{Ste} we have the following.

\begin{prop}\label{flat} $R_n$ is a flat left
$\bbz\G_n$--module so $\SA_n\cong\SA^\bbz_n\ox_{\bbz\G_n}R_n$.
Moreover $\SK_n$ is a flat $R_n$--module so
$\SA_n\ox_{R_n}\SK_n=H_1(M;\SK_n)$.
\end{prop}

Now we establish that the $R_n$ are PID's. Consider the short
exact sequence
$1\lra\wt G\lra G/\gs n_\bbq\xrightarrow{\pi}\bbz\lra1$ where
$\pi$ is induced by abelianization and $\wt G$ is the kernel
of $\pi$. Note that there are precisely two such epimorphisms
$\pi$. If we choose $\mu\in G/\gs n_\bbq$ which generates the
torsion-free part of the abelianization then $\pi$ is
canonical (take $\pi(\mu)=1$) and has a canonical splitting
$(1\xrightarrow{s}\mu)$. Now note that any element of
$\bbq[G/\gs n_\bbq]$ has a unique expression of the form
$\g=\mu^{-m}a_{-m}+\dots+a_0+\dots+\mu^ka_k$ where
$a_i\in\bbq\wt G$ ($a_{-m}$ and $a_k$ not zero unless
$\g=0$). Thus $\bbq[G/\gs n_\bbq]$ is canonically isomorphic
to the {\it skew Laurent polynomial ring}, $\bbq\wt
G[t^{\pm1}]$, in one variable with coefficients in $\bbq\wt
G$. Recall that the latter is the ring consisting of
expressions $t^{-m}a_{-m}+\dots+t^ka_k$, $a_i\in\bbq\wt G$
which add as ordinary polynomials but where multiplication is
twisted by an automorphism $\a\co\bbq\wt G\lra\bbq\wt G$ so
that if $a\in\bbq\wt G$ then $t^ia\cd t=t^{i+1}\a(a)$. The
automorphism in our case is induced by the automorphism of
$\wt G$ given by conjugation by $\mu$. The twisted
multiplication is evident in $\wt G$ since
$\mu^ia\cd\mu=\mu^i\mu(\mu^{-1}a\mu)=\mu^{i+1}\a(a)$.

Since $\wt G$ is a subgroup of a PTFA group, it also is PTFA
and so $\bbz\wt G$ admits a (right) skew field of fractions
$\bbk$ into which it embeds. This is also written $(\bbz\wt
G)(\bbz\wt G)^{-1}$ meaning that all the {\it non-zero}
elements of $\bbz\wt G$ are inverted. It follows that
$\bbz[G/\gs n_\bbq](\bbz\wt G)^{-1}$ is canonically identified
with the skew polynomial ring $\bbk[t^{\pm1}]$ with
coefficients in the skew field $\bbk$ (see
\cite[Proposition~3.2]{COT1} for more details). The following
is well known (see Chapter 3 of \cite{J} or Prop. 2.1.1 of
\cite{Co1}).

\begin{prop}\label{pid} A skew polynomial ring
$\bbk[t^{\pm1}]$ over a division ring $\bbk$ is a right (and
left) PID.
\end{prop}

\begin{proof} One first checks that there is a well-defined
degree function on any skew Laurent polynomial ring (over a
domain) where $\deg(t^{-m}a_{-m}+\dots+t^ka_k)=m+k$ and that this
degree function is additive under multiplication of polynomials.
Then one verifies that there is a division algorithm such that if
$\deg(q(t))\ge\deg(p(t))$ then $q(t)=p(t)s(t)+r(t)$ where
$\deg(r(t))<\deg(p(t))$. Finally, if $I$ is any non-zero right
ideal, choose $p\in I$ of minimal degree. For any $q\in I$,
$q=ps+r$ where, by minimality, $r=0$. Hence $I$ is principal. Thus
$\bbk[t^{\pm1}]$ is a right PID. The proof that it is a left PID
is identical.
\end{proof}

\begin{prop}\label{gamma-n} For $n\ge0$ let $R_n$ denote the
ring $\bbz[G/\gns+](\bbz\wt G)^{-1}$. This can be identified
with the PID $\bbk_n[t^{\pm1}]$ where $\bbk_n$ is the quotient
field of $\bbz\wt G$ ($1\lra\wt G\lra
G/\gns+\xrightarrow{\pi}\bbz\lra1$).
\end{prop}

Of course the isomorphism type of $\SA_n(K)$ is still purely a
function of the isomorphism type of the group $G$ of the knot
since $\SA_n(K)=\gns+/G^{(n+2)}\ox R_n$. However, when viewed as a
module over $\bbk_n[t^{\pm1}]$, it is also dependent on a choice
of the meridional element $\mu$. 

\medskip
\sh{Non-triviality}

We now show that the higher-order Alexander modules are {\it
never} trivial except when $K$ is a knot with Alexander
polynomial $1$. The following results generalize
Proposition~\ref{rank} and Lemma~\ref{2complex}.

\begin{cor}\label{rank estimate} If $X$ is a (possibly
infinite) $2$--complex with $H_2(X;\bbq)=0$ and
$\phi\co\pi_1(X)\lra\G$ is a PTFA coefficient system then
$\rank(H_1(X;\bbz\G))\ge\b_1(X)-1$.
\end{cor}

\begin{cor}\label{non-triviality} If $K$ is a knot whose
Alexander polynomial $\Delta_0$ is {\bf not} $1$, then the derived
series of $G=\pi_1(\sk)$ does not stabilize at finite $n$, i.e.\
$\gs n/\gns+\neq0$. Hence the derived module $\SA^\bbz_n(K)$ is
non-trivial for any $n$. Moreover, if $n>0$, $\SA_n(K)$ (viewed as
a $\bbk_n[t^{\pm1}]$ module) has rank at least
$\deg(\Delta_0(K))-1$ as a $\bbk_n$--module and hence is an
infinite dimensional $\bbq$ vector space.
\end{cor}

The first part of the Corollary has been independently established
by S.K. Roushon \cite {Ru}.

\medskip
{\bf Proposition~\ref{injectivity} $\Rightarrow$
Corollary~\ref{rank estimate}}\qua First consider the case that
$\b_1(X)$ is finite. Consider the case of
Proposition~\ref{injectivity} where $A$ is a wedge of
$\b_1(X)$ circles and $i\co A\lra X$ is chosen to be a
monomorphism on $H_1(\un{\ \ };\bbq)$. Then
$\rank(H_1(X;\bbz\G))$ is at least $\rank(H_1(A;\bbz\G))$
which is $\b_1(X)-1$ by Lemma~\ref{2complex}. Now if
$\b_1(X)$ is infinite, apply the above argument for a wedge
of $n$ circles where $n$ is arbitrary.

\medskip
{\bf Proposition~\ref{injectivity} $\Rightarrow$
Corollary~\ref{non-triviality}}\qua Let
$X$ be the infinite cyclic cover of $\sk$, and let $\wt
G=\pi_1(X)/\pi_1(X)^{(n)}=\gs1/\gns+$ as in
Proposition~\ref{gamma-n}. If $\Delta_0\neq1$ then
$\deg(\Delta_0)=\b_1(X)\ge2$. Applying Corollary~\ref{rank
estimate} we get that $H_1(X;\bbz\wt G)$ has rank at least
$\b_1(X)-1$. But $H_1(X;\bbz\wt G)$ can be interpreted as the
first homology of the $\wt G$--cover of $X$, as a $\bbz\wt G$
module. This covering space has $\pi_1$ equal to $\gns+$.
Since the $\wt G$ cover of $X$ is the same as the cover of
$\sk$ induced by $G\lra G/\gns+$, $H_1(X;\bbz\wt G)\cong
H_1(\sk;\bbz[G/\gns+])\equiv\SA^\bbz_n(K)$ as
$\bbz\wt G$--modules. Now, since $\SA^\bbz_n$ has rank at
least $\b_1(X)-1$ as a $\bbz\wt G$--module, $\SA_n$ has rank at
least $\b_1(X)-1$ as a $\bbk_n$ module since the latter is the
definition of the former. It follows that $\gns+/G^{(n+2)}$
is non-trivial (and hence infinite) for $n\ge0$. If $n>0$ it
follows that $\wt G$ is an infinite group. In this case
$\bbq\wt G$ and hence $\bbk_n$ are infinitely generated
vector spaces.

\section{Higher order Alexander
polynomials}\label{polynomials}

In this section we further analyze the localized Alexander modules
$\SA_n(K)$ that were defined in Section~\ref{rational} as right
modules over the skew Laurent polynomial rings
$R_n\cong\bbk_n[t^{\pm1}]$. We define higher-order ``Alexander
polynomials'' $\Delta_n(K)$ and show that their degrees $\d_n(K)$
are integral invariants of the knot. We prove that $\d_0$,
$\d_1+1$, $\d_2+1,\dots$ is a non-decreasing sequence for any
knot. In later sections we will see that the $\d_n$ are powerful
knot invariants with applications to genus and fibering questions.
The higher-order Alexander polynomials bear further study.

Recall that it has already been established that $\SA_n(K)$
is a finitely-generated torsion right $R_n$ module where
$R_n$ is a PID. The following generalization of the standard
theorem for commutative PIDs is well known (see Theorem 2.4
p. 494 of \cite{Co2}).

\begin{thm}\label{classification} Let $R$ be a principal
ideal domain. Then any finitely generated torsion right
$R$--module $M$ is a direct sum of cyclic modules
$$
M\cong R/e_1R\op\dots\op R/e_rR
$$
where $e_i$ is a total divisor of $e_{i+1}$ and this
condition determines the $e_i$ up to similarity.
\end{thm}

Here $a$ is {\it similar} to $b$ if $R/aR\cong R/bR$ (p. 27
\cite{Co1}). For the definition of {\it total divisor}, the
reader is referred to Chapter 8 of \cite{Co2}. This
complication is usually unnecessary because a finitely
generated torsion module over a {\it simple} PID is cyclic
(pp. 495--496 \cite{Co2})!! For $n>0$, $R_n$ is almost always
a simple ring, but since this fact will not be used in this
paper, we do not justify it.

\begin{defn}\label{higher-order polynomials} For any knot $K$
and any integer $n\ge0$, $\{e_1(K),\dots,e_r(K)\}$ are the
elements of the PID $R_n$, well-defined up to similarity,
associated to the canonical decomposition of $\SA_n(K)$. Let
$\Delta_n(K)$, the {\it $n^\supth$ order Alexander polynomial} of
$K$, be the product of these elements, viewed as an element of
$\bbk_n[t^{\pm1}]$ (for $n=0$ this is the classical Alexander
polynomial).
\end{defn}

The polynomial $\Delta_n(K)$, as an element of $R_n$, is also
well-defined up to similarity (a non-obvious fact that we will not
use). However as an element of $\bbk_n[t^{\pm1}]$ it acquires
additional ambiguity because a splitting of $G\thra\bbz$ was used
to choose an isomorphism between $R_n$ and $\bbk_n[t^{\pm1}]$.
Alternatively, using a square presentation matrix for $\SA_n(K)$
(see the next section), one can associate an element of $K_1(R_n)$
 and, using the Dieudonn\'{e} determinant, recover $\Delta_n(K)$ as an element of $U/[U,U]$ where $U$
  is the group of units of the quotient field of $R_n$. Since similarity is not well-understood in a noncommutative ring
(being much more difficult than merely identifying when elements
differ by units), we have not yet been able to make effective use
of the higher-order Alexander polynomials except for their
degrees, which turn out to be perfectly well-defined integral
invariants, as we now explain.

\begin{defn}\label{delta-n} For any knot $K$ and any integer
$n\ge0$, the {\it degree} of the $n^\supth$ order Alexander
polynomial, denoted $\d_n(K)$ is an invariant of $K$. It can
be defined in any of the following equivalent ways:
\begin{enumerate}
\item[1)] the degree of $\Delta_n(K)$
\item[2)] the sum of the degrees of $e_i(K)\in
R_n\cong\bbk_n[t^{\pm1}]$
\item[3)] the rank of $\SA_n(K)$ as a module over $\bbk_n$
\item[4)] the rank of $\gns+/G^{(n+2)}\ox_{\bbz\G_n}R_n$ as a
module over the subring $\bbz\wt G\sbq\bbz\G_n$
\item[5)] the rank of $\gns+/G^{(n+2)}$ as a module over the
subring $\bbz[\gs1/\gns+]\subset\bbz[G/\gns+]$.
\end{enumerate}
\end{defn}

\begin{proof}[Proof of Definition~\ref{delta-n}] Definitions 4
and 5 are independent of choices since there $R_n$ has not been
specifically identified with the polynomial ring $\bbk_n$. To see
that Definition 3 is the same as 4, consider
Definition~\ref{module} and Proposition~\ref{flat}. Also note that
the identification of $\bbz[G/\gns+]$ with the skew polynomial
ring $\bbz\wt G[t^{\pm1}]$, carries the subring $\bbz\wt G$ ({\it
independent of splitting}) to the ring of elements of degree zero.
Thus under any identification of $R_n=\bbz[G/\gns+]$ $(\bbz\wt
G-\{0\})^{-1}$ with $\bbk_n[t^{\pm1}]$, the quotient field
$\bbz\wt G(\bbz\wt G-\{0\})^{-1}$ is carried (independent of
splitting) to $\bbk_n$, viewed as the subfield of elements of
degree zero. From Definition 3 and Theorem~\ref{classification},
one sees that these ranks are finite because the rank of
$\bbk_n[t^{\pm1}]/p(t)\bbk_n[t^{\pm1}]$ is easily seen to be the
degree of $p(t)$. The equivalence of Definitions 1 and 2 then
follows trivially. To see that 4 and 5 are equivalent, one must
show that $\SA^\bbz_n\ox_{\bbz[G/\gns+]}\bbk_n[t^{\pm1}]$ as a
$\bbk_n$--module is merely $\SA^\bbz_n\ox_{\bbz\wt G}\bbk_n$. This
is left to the reader.
\end{proof}

We can establish one interesting property of the $\d_n$,
namely that for any $K$ they form a non-decreasing sequence.
This theorem says that the derived series of the fundamental
group of a knot complement (more generally of certain
$2$--complexes) cannot stabilize unless $\d_0=1$ (see
Corollary~\ref{non-triviality}). Moreover in some sense the
``size'' of the successive quotients $G^{(n)}/\gns+$ is
non-decreasing.

\begin{thm}\label{non-decreasing} If $K$ is a knot then
$\d_0(K)\le\d_1(K)+1\le\d_2(K)+1\le\dots\le\d_n(K)+1$.
\end{thm}

\begin{proof} First we show $\d_1\ge\d_0-1$. Let $X$ be the
infinite cyclic cover of $\sk$ and $G=\pi_1(\sk)$. Note that
$\b_1(X)=\rank_\bbq H_1(\sk;\bbq[t,t^{-1}])=\d_0$, and
$\d_1=\rank_{\bbk_1}H_1(\sk;\bbk_1[t^{\pm1}])=\rank
H_1(X;\bbz[\gs1/\gs2])$ by Definition~\ref{delta-n}. The latter,
by Corollary~\ref{rank estimate} is at least $\b_1(X)-1$ (since
$H_2(X;\bbq)=0)$ and we are done.

Now it will suffice to show $\d_n\ge\d_{n-1}$ if $n\ge2$. Let
$X_n$ be the covering space of $\sk$ with fundamental group
$\gns+$ so $X_0=X$. Then $X_{n-1}$ is a covering space of $X$
with $\gs1/\gs n$ as deck translations. Choose a wedge of
$\d_0$ circles $A_0\to X$ giving an isomorphism on $H_1(\un{\
\ };\bbq)$. Let $\wt A_0\xrightarrow{i}X_{n-1}$ be the
induced cover and corresponding inclusion. By
Proposition~\ref{injectivity}, $i_*$ is a monomorphism on
$H_1$. Since $n\ge2$,
$\rank_{\bbz[\gs1/\gs n]}H_1(A_0;\bbz[\gs1/\gs n])$ is
precisely $\b_1(A_0)-1=\d_0-1$ by Lemma~\ref{2complex} (here
we assume $\d_0>0$ since if $\d_0=0$ then $\d_i=0$ and the
theorem holds). Choose a subset of image $i_*$ with
cardinality $\d_0-1$ that is $\bbz[\gs1/\gs n]$--linearly
independent in $H_1(X_{n-1})$. It is not difficult to show
that, in a module over an Ore domain, any linearly
independent set can be extended to a maximal linearly
independent set, i.e.\ whose cardinality is equal to the rank
of the module. Hence if $\d_{n-1}$ (which equals the
$\bbz[\gs1/\gs n]$--rank of $H_1(X_{n-1})$) exceeds $\d_0-1$,
then there is a set of $e=\d_{n-1}-(\d_0-1)$ circles and a map
$\tl f\co A_e\to X_{n-1}$ of a wedge of $e$ circles, such that
the free submodule generated by these circles captures the
``excess rank.'' Let $f=\pi\circ\tl f\co A_e\to X$. Then the map
$A=A_0\vee A_e\lra X$ induces a monomorphism on
$H_1(\un{\ \ };\bbz[\gs1/\gs n])$ by construction. Another
way of saying this is that the induced map on $\gs1/\gs
n$--covers $A_{n-1}\to X_{n-1}$ is injective on $H_1(\un{\ \
};\bbz)$ where $A_{n-1}$ is the induced cover of $A$. Since
$H_2(X;\bbz)=0$, it follows from Lemma~\ref{2complex} that
$H_2(X_{n-1};\bbz)=0$. Hence $(X_{n-1},A_{n-1})$ is a
relative $2$--complex that satisfies the conditions of
Proposition~\ref{injectivity}, with $\G=\gs n/\gns+$. It
follows that
$H_1(A_{n-1};\bbz\G)\xrightarrow{i_*}H_1(X_{n-1};\bbz\G)$ is
injective. But this is the same as the map induced by $i\co A\to
X$ on $H_1(\un{\ \ };\bbz[\gs1/\gns+])$. Thus $\d_n=\rank
H_1(X;\bbz[\gs1/\gns+])$ is at least the rank of
$H_1(A;\bbz[\gs1/\gns+])$. Since $A$ is a wedge of
$e+\d_0=\d_{n-1}+1$ circles and $n\ge2$, this latter rank is
precisely $\d_{n-1}$ by Lemma~\ref{2complex}. Hence
$\d_n\ge\d_{n-1}$ as claimed.
\end{proof}

\begin{ques} Is there a knot $K$ and some $n>0$ for which
$\d_n(K)$ is a non-zero even integer?
\end{ques}

If not then a complete realization theorem for the $\d_i$ can
be derived from the techniques of Section~\ref{examples}.

\section{Presentation of $\SA_n$ from a Seifert
surface}\label{seifert}

Suppose $M$ is a knot exterior, or more generally a compact,
connected, oriented $3$--manifold with $\b_1=1$ that is either
closed or whose boundary is a torus. Suppose $\pm V$ is a compact,
connected, oriented surface which generates $H_2(M,\p M)$. In the
case of a knot exterior, the orientation on the knot can be used
to fix the orientation of $V$, and $V$ can be chosen to be a {\it
Seifert surface\/} of $K$. The classical Alexander module of $K$
can be calculated from a presentation matrix which is obtained by
pushing certain loops in $V$ into $S^3\setminus V$. Here we show
that there is a finite presentation of $\SA_n(K)$ obtained in a
similar fashion from $V$.

Let $Y=M-(V\x(-1,1))$ and denote by $i_+$ and $i_-$ the two
inclusions $V\lra V\x\{\pm1\}\lra\p Y\subset Y$. Recall from
Definition~\ref{an} and Proposition~\ref{gamma-n} that
$\SA_n(M)\cong H_1(M;\bbk_n[t^{\pm1}])$ where an isomorphism is
fixed by choosing a circle $u$ dual to $V$ (an oriented meridian
in the case that $M=\sk$). The derivation of a presentation for
$\SA_n(M)$ follows the classical case (see page~122--123 of
\cite{Hi2})but is complicated by basepoint concerns. The following
overlaps with work of S. Harvey \cite{Ha}.

\begin{prop}\label{presentation} The following sequence is exact.
$$
H_1(V;\bbk_n)\ox_{\bbk_n}\bbk_n[t^{\pm1}]\xrightarrow{d}
H_1(Y;\bbk_n)\ox_{\bbk_n}\bbk_n[t^{\pm1}]\lra\SA_n(M)\lra0
$$
where $d(\a\ox1)=(i_+)_*\a\ox t-(i_-)_*\a\ox1$.
\end{prop}

\begin{proof} (see {\rm\cite{Ha}} for a more detailed proof) For simplicity let $\G_n$ stand for
$G/\gns+_{\bbq}$ so there is an exact sequence $1\lra\wt
G\lra\G_n\xrightarrow{\pi}\bbz\lra1$ where $\pi(u)=1$ and $\bbk_n$
is the quotient (skew) field of $\bbz\wt G$. Let $U=V\x[-1,1]$ and
consider a Mayer--Vietoris sequence for homology with $\bbz\G_n$
coefficients using the decomposition $M=Y\cup U$. Or, more
naively, consider an ordinary Mayer--Vietoris sequence for the
integral homology of $M_{\G_n}$, the $\G_n$ cover, using the
decomposition $M_{\G_n}=p^{-1}(Y)\cup p^{-1}(U)=Y_{\G_n}\cup
U_{\G_n}$ and note that all the maps are $\bbz\G_n$--module
homomorphisms. After the usual simplification one arrives at the
exact sequence:
$$
\lra H_1(V;\bbz\G_n)\xrightarrow{d}H_1(Y;\bbz\G_n)
\xrightarrow{j_*}\SA^\bbz_n(M)\xrightarrow{\p_*}H_0(V;\bbz\G_n).
$$
Localizing yields a similar sequence with $\bbk_n[t^{\pm1}]$
coefficients where $\SA_n(M)$ replaces $\SA^\bbz_n(M)$. Since
$\pi_1(V)$ and $\pi_1(Y)$ are contained in $\wt G$, one can
consider $H_*(V;\bbk_n)$ and $H_*(Y;\bbk_n)$, which are free
$\bbk_n$--modules. Moreover $\bbk_n[t^{\pm1}]$ is free and
hence flat as a left $\bbk_n$ module. Thus
$H_*(V;\bbk_n[t^{\pm1}])\cong
H_*(V;\bbk_n)\ox_{\bbk_n}\bbk_n[t^{\pm1}]$ and
$H_*(Y;\bbk_n[t^{\pm1}])\cong
H_*(Y;\bbk_n)\ox_{\bbk_n}\bbk_n[t^{\pm1}]$, showing that
these homology groups are finitely-generated {\it free}
$\bbk_n[t^{\pm1}]$ modules. Since $\SA_n(M)$ is a torsion
module by Proposition~\ref{finite generation} and
$H_0(V;\bbk_n[t^{\pm1}])$ is free, $\p_*$ is the zero map.
This concludes our sketch of the proof of the proposition.
\end{proof}

\begin{cor}\label{square matrix} If the (classical) Alexander polynomial of
$M$ is not $1$ then $\SA_n(M)$, $n>0$, has a square presentation matrix of
size $r=\max\{0,-\chi(V)\}$ each entry of which is a Laurent polynomial of
degree at most $1$. Specifically, we have the presentation
$$
(\bbk_n[t^{\pm1}])^r\xrightarrow{\p}
(\bbk_n[t^{\pm1}])^r\lra\SA_n(M)\lra0
$$
where $\p$ arises from the above proposition. If $n=0$ then
the same holds with $r$ replaced by $\b_1(V)$.
\end{cor}

\begin{proof} The Corollary will follow immediately from the
Proposition if we establish that $H_1(V;\bbk_n)\cong
H_1(Y;\bbk_n)\cong\bbk^r_n$. Note that both $V$ and $Y$ have
the homotopy type of finite connected $2$--complexes. Consider
the coefficient systems $\psi\co\pi_1(V)\lra\wt G$ and
$\psi'\co\pi_1(Y)\lra\wt G$ obtained by restriction of
$\pi_1(M)\lra\G_n$. Letting $b_i$ stand for the rank of
$H_i(\un{\ \ };\bbz\wt G)$ or equivalently the rank of
$H_i(\un{\ \ };\bbk_n)$, we have that
$\chi(V)=b_0(V)-b_1(V)+b_2(V)$ as in Fact 3.

Suppose that $\psi$ is non-trivial. Then $b_0(V)=0$ by
Proposition~\ref{Ho}. Since $\wt G$ is PTFA, it is torsion free
and hence the image of $\psi$ is infinite. It follows that the
$\wt G$--cover of $V$ is a non-compact $2$--manifold and thus
$b_2(V)=0$. Therefore $b_1(V)=r$ as desired. It also follows that
$\psi'$ is non-trivial and so $b_0(Y)=0$. Since $\chi(M)=0$ it
follows that $\chi(Y)=\chi(V)$. Thus
$b_2(Y)-b_1(Y)=\chi(Y)=\chi(V)=-b_1(V)$ so $b_1(Y)=b_1(V)+b_2(Y)$.
By Proposition~\ref{presentation} $\SA_n$ has a presentation of
deficiency $b_1(Y)-b_1(V)$. If $b_2(Y)>0$ then $\SA_n(M)$ has a
presentation of positive deficiency, contradicting the fact that
it is a $\bbk_n[t^{\pm1}]$--torsion module. Therefore $b_2(Y)=0$
and $b_1(Y)=b_1(V)=r$ as required. This completes the case that
$\psi$ is non-trivial, after noting that if $n=0$ then $\psi$ is
certainly trivial since $\wt G=1$.

Now suppose $\psi$ is trivial. If $n=0$ then this is the case
of the classical (rational) Alexander module and the result
is well-known. If $n\ge1$ then the triviality of $\psi$
implies that $\pi_1(V)\sbq\gs2_\bbq$. Consider a map $f\co M\lra
S^1$ such that $V$ is the inverse image of a regular value.
Then $\gs1_\bbq=\ker f_*$ and it follows that $\gs1_\bbq$ is
the normal subgroup generated by $\pi_1(Y)$ and so, for any
$\g\in\pi_1(Y)$, there exists a non-zero integer $m$ such
that $m\g$ bounds an orientable surface $S$. Hence
$\gs1_\bbq/\gs2_\bbq$ is generated by $\pi_1(V)$ and thus is
zero. It follows that $\SA_0(M)=0$ and that classical
Alexander polynomial is $1$. Since this case was excluded by
hypothesis, the proof is complete.
\end{proof}

\begin{ex}\label{fibered knot} Suppose $K$ is a fibered knot of genus
$g$ with fiber surface $V$ and $\pi_1$--monodromy $f$. If $n>0$ and
 $F$ is the free group of rank $2g-1$ then $H_1(V;\bbk_n)\cong
H_1(Y;\bbk_n)\cong H_1(F;\bbk_n)\cong\bbk_n^{2g-1}$ by
Lemma~\ref{2complex}. By the above results, $\SA_n$ has a $(2g-1)$
by $(2g-1)$ presentation matrix given by $It-f_n$ where $f_n$ is
an automorphism of the vector space $\bbk_n^{2g-1}$ derived from
the induced action of $f$ on $F/F^{(n+1)}$ .
\end{ex}

\section{The $\d_n$ give lower bounds for knot
genus}\label{examples}

The previous section can now be used to show that the degrees of
the higher order Alexander polynomials give lower bounds for
$\genus(K)$. In the last part of this section we show that there
are knots such that $\d_0<\d_n+1$ so that these invariants yield
sharper estimates of knot genus than that given by the Alexander
polynomial, $\deg(\Delta_0)\le2\genus(K)$. S. Harvey has
established analagous results for any $3$--manifold, finding
lower-bounds for the Thurston norm \cite{Ha}.

\begin{thm}\label{genus} If $K$ is a (null-homologous)
non-trivial knot in a rational homology sphere and $\d_n$ is
the degree of the $n^{\supth}$ order Alexander polynomial
then $\d_0\le2\genus(K)$ and $\d_n+1\le2\genus(K)$ if $n>0$.
\end{thm}

\noindent{\it Proof}. We may assume $n>0$ since the result for
$n=0$ is well known. If the classical Alexander polynomial is $1$
then $\d_0=\d_n=0$ and the theorem holds. Otherwise suppose $V$ is
a Seifert surface of minimal genus. By Corollary~\ref{square
matrix} $\SA_n(K)$ has a square presentation matrix of size
$2\genus(K)-1$. Since $\d_n$ is defined as $\rank_{\bbk_n}\SA_n$,
it remains only to show that the latter is at most $2\genus(K)-1$.
This is accomplished by the following lemma of Harvey.

\begin{lem}\label{determinant} {\rm\cite{Ha}}\qua Suppose $\SA$ is
a torsion module over a skew Laurent polynomial ring
$\bbk[t^{\pm1}]$ where $\bbk$ is a division ring. If $\SA$ is
presented by an $m\x m$ matrix $\th$ each of whose entries is
of the form $ta+b$ with $a$, $b\in\bbk$, then the rank of
$\SA$ as a $\bbk$--module is at most $m$.
\end{lem}


\begin{thm}\label{increasing degree} For any knot $K$ whose
(classical) Alexander polynomial is not $1$ and any positive
integer $k$, there exists a knot $K_*$ such that
\begin{description}
\item[a)]  $\SA_n(K_*)\cong\SA_n(K)$ for all $n<k$.
\item[b)] $\d_n(K_*)=\d_n(K)$ for all $n<k$.
\item[c)] $\d_k(K_*)>\d_k(K)$.
\item[d)] $K_*$ can be taken to be hyperbolic or to be
concordant to $K$.
\end{description}
\end{thm}

\begin{cor}\label{increasing degrees} Under the hypotheses of
the theorem above, there exists a hyperbolic knot $K_*$, with
the same classical Alexander module as $K$, for which
$\d_0(K_*)<\d_1(K_*)<\dots<\d_k(K_*)$.
\end{cor}

\begin{proof}[Proof of Theorem~\ref{increasing degree}] Let
$P=\pi_1(\sk)$ and let $\a$ be an element of
$P^{(k)}$ which does not lie in $P^{(k+1)}$. By
Corollary~\ref{non-triviality} such $\a$ are plentiful. We now
describe how to construct a knot $K_*=K(\a,k)$ which differs from
$K$ by a single ``ribbon move,'' i.e.\ $K_*$ is obtained by
adjoining a trivial circle $J$ to $K$ and then fusing $K$ to this
circle by a band as shown in Figure~\ref{example}. Thus $K_*$
is concordant to $K$. From a group theory perspective, what
is going on is simple. It is possible to add one generator
and one relation that precisely kills that generator if one
``looks'' modulo $n^{\supth}$ order commutators, but does not
kill that generator if one ``looks'' modulo $(n+1)^{\rm st}$ order
commutators. Details follow.

\begin{figure}[ht!]
\cl{\relabelbox\small
\epsfxsize 4in \epsfbox{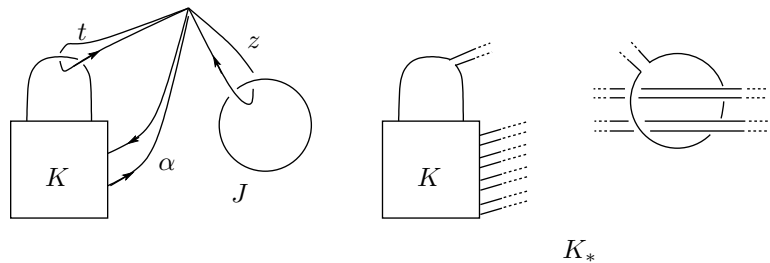}
\relabel {K1}{$K$}
\relabel {K2}{$K$}
\relabel {K3}{$K_*$}
\relabel {J}{$J$}
\relabel {t}{$t$}
\relabel {z}{$z$}
\relabel {a}{$\alpha$}
\endrelabelbox}
\caption{$K_*$ is obtained from $K$ by a ribbon move}
\label{example}
\end{figure}

Choose meridians $t$ and $z$ as shown. Choose an embedded band
which follows an arc in the homotopy class of the word
$\eta=t[\a^{-1},t^{-1}z]t^{-1}$. There are many such bands.
For simplicity choose one which pierces the disk bounded by
$J$ precisely twice corresponding to the occurrences of $z$
and $z^{-1}$ in $\eta$. Let $G=\pi_1(S^3-K_*)$ and let $\g$
denote a small circle which links the band. A Seifert
Van--Kampen argument yields that the group $E\equiv G/\<\g\>$
has a presentation obtained from a presentation of $P$ by
adding a single generator $z$ (corresponding to the meridian
of the trivial component) and a single relation $z=\eta
t\eta^{-1}$. We symbolize this by
$E=\<P,z\mid z=\eta t\eta^{-1}\>$. First we analyze the
relationship between $P$ and $E$.

\begin{lem}\label{increasing deg} Given $P$, $\a$, $k$,
$t$, $z$, $E$ as above:

{\rm a)}\qua $P/P^{(n)}\cong E/E^{(n)}$ for all $n\le k+1$ implying
that for all $n<k$, $\SA^\bbz_n(P)\cong\SA^\bbz_n(E)$ and
$\d_n(P)=\d_n(E)$;

{\rm b)}\qua $\d_k(E)=\d_k(P)+1$.
\end{lem}

\begin{proof}[Proof of Lemma~\ref{increasing deg}] Let
$w=t^{-1}z$ so $E=\<P,w\mid w=[t^{-1},\eta]\>$ and
$\eta=t[\a^{-1},w]t^{-1}$. Since $\a\in P^{(k)}$, $\eta\in
E^{(k)}$ and hence $w\in E^{(k)}$. But then $\eta\in
E^{(k+1)}$ so $w\in E^{(k+1)}$. Part a) of the Lemma follows
immediately: the epimorphism $E\lra P$ obtained by killing $w$
induces an isomorphism $E/E^{(k+1)}\lra P/P^{(k+1)}$, and
hence $\SA^\bbz_n(P)\cong\SA^\bbz_n(E)$ for $n<k$ by
Definition~\ref{general def}. Here we use the fact that both
$E$ and $P$ are {\it$\SE$--groups} in the sense of R.~Strebel
(being fundamental groups of $2$--complexes with $H_2=0$ and
$H_1$ torsion-free). Consequently any term of their derived
series is also an $\SE$--group and it follows that their
derived series is identical to their rational derived series
\cite{Str}.

Now we consider the subgroup $E^{(k+1)}$ of $E$. To justify
the following group-theoretic statements, consider a
$2$--complex $X$ whose fundamental group is $P$ and define a
$2$--complex $Y$ by adjoining a $1$--cell and a $2$--cell so
that $\pi_1(Y)\cong E$ corresponding to the presentation
$\<P,w\mid w=[t^{-1},\eta]\>$. The subgroup $E^{(k+1)}$ is
thus obtained by taking the infinite cyclic cover $Y_\infty$
of $Y$ (so $\pi_1(Y_\infty)=E^{(1)}$) followed by taking the
$E^{(1)}/E^{(k+1)}$--cover $\wt Y$ of $Y_\infty$ (so
$\pi_1(\wt Y)=E^{(k+1)}$). Since the inclusion map $X\lra Y$
induces an isomorphism $P/P^{(k+1)}\lra E/E^{(k+1)}$, the
induced cover of the subspace $X\subset Y$ is the cover $\wt
X$ of $X$ with $\pi_1(X)\cong P^{(k+1)}$. Therefore a cell
structure for $\wt Y$ relative to $\wt X$ contains only the
lifts of the $1$--cell $w$ and the $2$--cell corresponding to
the single relation. This allows for an elementary analysis
of $E^{(k+1)}$ as follows. By analyzing $X_\infty$ and
$Y_\infty$ we see that
$$
E^{(1)} = \<P^{(1)},\ w_i\ i\in\bbz\mid w_i = t^{-i}
[t^{-1},\eta] t^i\>
$$
where $w_i$ stands for $t^{-i}wt^i$ as an element of
$\pi_1(Y)$. If we rewrite the relation using
$\b^{-1}=t^{-i}\a^{-1}t^i$ and
$r^{-1}=t^{-i+1}\a^{-1}t^{i-1}$ we get
$$
E^{(1)} = \<P^{(1)},\ w_i\mid w_i = \b^{-1}w_i\b w^{-1}_i
w_{i-1} r^{-1}w^{-1}_{i-1} r\>.
$$
This is a convenient form because what we want to do now is
``forget the $t$ action'' because $\d_k$ is defined as the rank of
the abelianization of $E^{(k+1)}$ {\it as a module over}
$\bbz[E^{(1)}/E^{(k+1)}]$ (or equivalently over its quotient field
$\bbk_k$). Therefore we now think of $Y_\infty$ as being obtained
from $X_\infty$ by adding an infinite number of $1$--cells $w_i$
and a correspondingly infinite number of $2$--cells. Thus $\wt Y$
is obtained from $\wt X$ by adding $1$--cells $\{w^s_i\mid
i\in\bbz, s\in E^{(1)}/E^{(k+1)}\}$, where $w^s_i$ descends to
$s^{-1}t^{-i}wt^is$ in $E$, and $2$--cells corresponding to the
relations $\{w^s_i=w^{\b s}_i(w^s_i)^{-1}
w^s_{i-1}(w^{rs}_{i-1})^{-1}\mid i\in\bbz, s\in
E^{(1)}/E^{(k+1)}\}$ where, for example, $w^{\b s}_i$ is the image
of a fixed $1$--cell $w_i$ under the deck translation $\b s\in
E^{(1)}/E^{(k+1)}$ and descends to $s^{-1}\b^{-1}t^{-i}wt^i\b s$
in $E$. The abelianization, $E^{(k+1)}/E^{(k+2)}$, as a right
$\bbz[E^{(1)}/E^{(k+1)}]\cong\bbz[P^{(1)}/P^{(k+1)}]$ module is
obtained from $P^{(k+1)}/P^{(k+2)}$ by adjoining a generator $w_i$
and a relation for each $i\in\bbz$. Upon rewriting the relations
above as $w_i(2s-\b s)_*=w_{i-1}(s-rs)_*$ where $(2s-\b s)_*$
denotes the (right) action of $2s-\b s\in\bbz[P^{(1)}/P^{(k+1)}]$,
and then again as $w_i(2-\b)_*s_*=w_{i-1}(1-r)_*s_*$ we see that
the relations are generated as a module by
$\{w_i(2-\b)_*=w_{i-1}(1-r)_*\mid i\in\bbz\}$. Note that neither
$2-\b$ nor $1-r$ is zero since their augmentations are non-zero.
Hence in $\bbk_k$ these elements are invertible and each $w_i$,
$i\neq0$ can be equated uniquely to a multiple of $w_0$. Thus
$E^{(k+1)}/E^{(k+2)}\cong P^{(k+1)}/P^{(k+2)}\op\bbk_k$ as a
$\bbk_k$--module. It follows immediately that $\d_k(E)=\d_k(P)+1$.
This concludes the proof of Lemma~\ref{increasing deg}.
\end{proof}

Returning to the proof of the theorem, it will suffice to
show $\g\in G^{(k+1)}$ since if so then the epimorphism
$G\lra E$ induces an isomorphism $G/\gs n\cong E/E^{(n)}$ for
all $n\le k+1$ and hence an isomorphism
$\SA^\bbz_n(G)\lra\SA^\bbz_n(E)$ for $n<k$. Moreover the
epimorphism $G^{(k+1)}\lra E^{(k+1)}$ induces an epimorphism
$\SA^\bbz_k(G)\lra\SA^\bbz_k(E)$ of $G/G^{(k+1)}$ ($\cong
E/E^{(k+1)}$) modules. Thus $\d_k(K_*)=\d_k(G)\ge\d_k(E)$. By
Lemma~\ref{increasing deg} the map $P\lra E$ induces
isomorphisms $\SA^\bbz_n(P)\lra\SA^\bbz_n(E)$ for $n<k$ and
$\d_k(E)=\d_k(P)+1=\d_k(K)+1$. Combining these results will
finish the proof.

\begin{figure}[ht!]
\cl{\relabelbox\small
\epsfxsize 4in \epsfbox{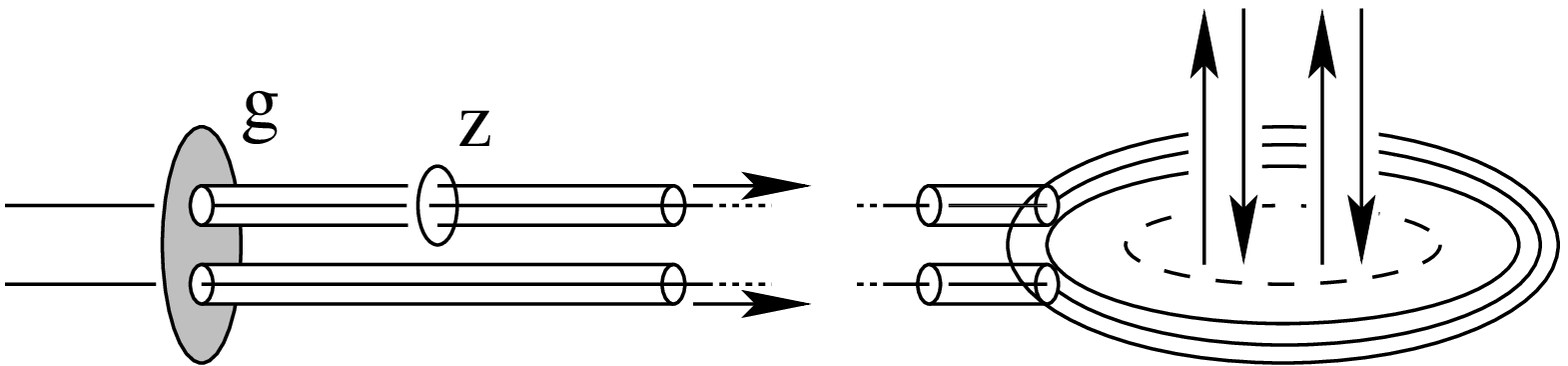}
\adjustrelabel <-5pt, 0pt> {g}{$\gamma$}
\adjustrelabel <-5pt, 0pt> {z}{$z^*$}
\extralabel <-1.3in, 0.28in> {$\ell_z$}
\endrelabelbox}
\caption{$\g=[z^*,\ell_z]$} \label{torus1}
\end{figure}

\begin{figure}[ht!]
\cl{\relabelbox\small
\epsfxsize 1.7in \epsfbox{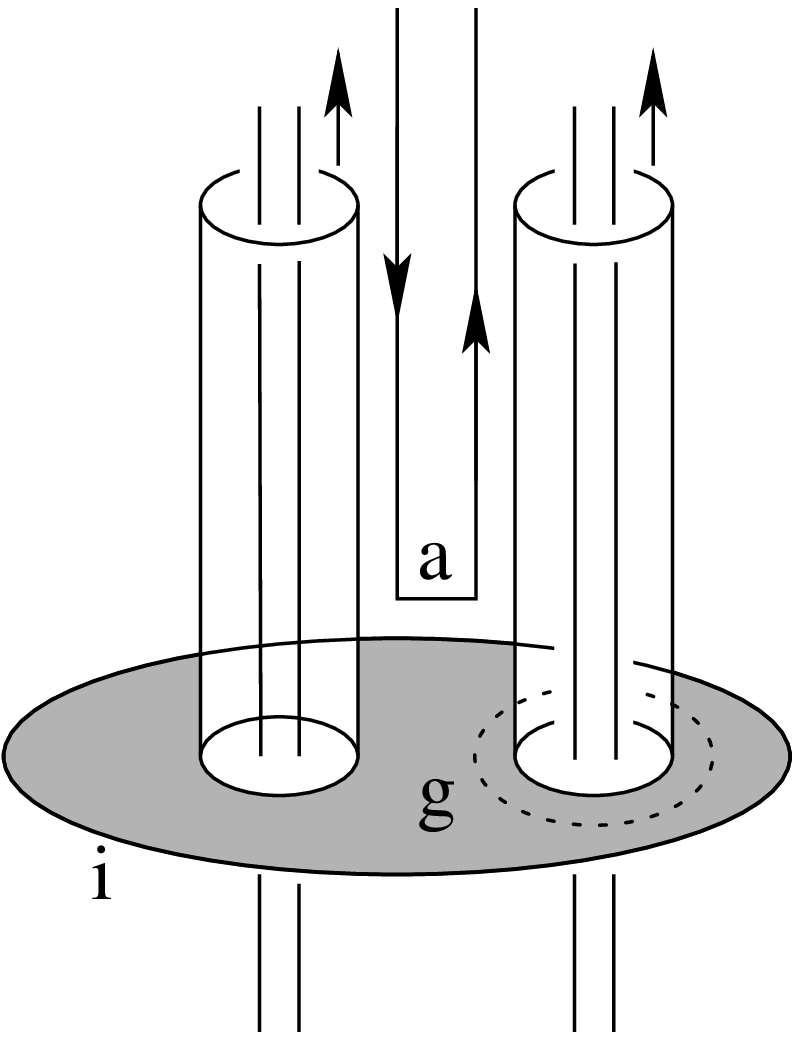}
\adjustrelabel <0pt, 0pt> {a}{$\bar\alpha$}
\adjustrelabel <-2pt, 0pt> {g}{$\gamma^*$}
\adjustrelabel <0pt, 0pt> {i}{$\ell_z$}
\endrelabelbox}
\caption{$\ell_z=[\g^*,\bar\a]$} \label{torus2}
\end{figure}

To see that $\g\in G^{(k+1)}$, first note that $\g$ bounds an
embedded disk which is punctured twice by the knot. By tubing
along the knot in the direction of $J$, one sees that $\g$
bounds an embedded (punctured) torus in $\sk_*$ as in
Figure~\ref{torus1}. This illustrates the group-theoretic fact
that $\g=[z^*,\ell_z]$ where $\ell_z$ is a longitude of $J$
and $z^*$ is a conjugate of $z$. It suffices to show
$\ell_z\in G^{(k+1)}$. But, since $\eta=t\a^{-1}t^{-1}z\a
z^{-1}$ contains 2 occurrences of $z$ (with opposite sign)
and we chose our band to pass precisely 2 times through $J$,
$\ell_z$ bounds a twice punctured disk and hence a punctured
torus as in Figure~\ref{torus2}. This illustrates that
$\ell_z=[\g^*,\bar\a]$ where $\g^*$ is a conjugate of
$\g$ since it is another meridian of the band, and $\bar\a$
is the word $\a$ separating the occurrences of $z$ and
$z^{-1}$ in the word $\eta$. Clearly $\g\in\gs1$. Suppose
$\g$, and hence $\g^*$, lies in $\gs j$ for some $1\le j\le
k$. Thus $G/\gs j\cong E/E^{(j)}$. Let $\a'$ denote the image
of $\bar\a$ under the map $G\to E$. Then $\a'$ is the image of
$\a$ under the map $P\to E$ since all the elements $\a$,
$\a'$ and $\bar\a$ are represented by the ``same'' path. Since
$\a\in P^{(k)}$ (by hypothesis), $\a'\in E^{(k)}$ and hence
$\bar\a\in\gs j$. But then $\ell_z\in G^{(j+1)}$ and hence
$\g\in G^{(j+1)}$. Continuing in this way shows that $\g\in
G^{(k+1)}$ and concludes the proof of Theorem~\ref{increasing
degree}.
\end{proof}

\begin{proof}[Proof of Corollary~\ref{increasing degrees}]
By induction and Theorem~\ref{increasing degree} there exists a
knot $K_{k-1}$ with the same classical Alexander module as $K$ and
$\d_0(K_{k-1})<...<\d_{k-1}(K_{k-1})$  Apply
Theorem~\ref{increasing degree} to $K_{k-1}$  produce a new knot
$K_*$. One easily checks that $K_*$ satisfies the required
properties by Theorem~\ref{increasing degree},
Theorem~\ref{non-decreasing} and Corollary~\ref{hyperbolic}.
\end{proof}

\section{Genetic infection: A technique for constructing
knots}\label{genetic}

We discuss a satellite construction, which we call {\it genetic
modification} or {\it infection}, by which a given knot $K$ is
subtly modified, or {\it infected} using an auxiliary knot or link
$J$ (see also  of \cite[Section 6]{COT1} \cite{COT2} \cite{CT}).
If, by analogy, we think of the {\it group} $G$ of $K$ as its {\it
strand of DNA}, then, by Corollary~\ref{non-triviality}, this
``strand'' is infinitely long as measured by the derived series.
Thus, as we shall see, it is possible to locate a spot on the
``strand'' which corresponds to an element of $\gs n-\gns+$,
excise a ``small piece of DNA'' and replace it with ``DNA
associated to the knot $J$'', with the effect that $G/\gns+$ is
not altered but $G/G^{(n+2)}$ is changed in a predictable fashion.
The infection is subtle enough so that it is not detected by the
localized modules $\SA_n$ (hence not by $\d_n$). The effect on the
(integral) modules $\SA^\bbz_n$ can be measured {\it numerically}
by the higher-order signatures of Section~\ref{signatures}.

Suppose $K$ and $J$ are fixed knots and $\eta$ is an embedded
oriented circle in $\sk$ which is itself unknotted in $S^3$. Note
that any class $[\eta]\in G$ has a (non-unique) representative
$\eta$ which is unknotted in $S^3$. Then $(K,\eta)$ is isotopic to
part a of Figure~\ref{infection} below, where some
undetermined number $m$ of strands of $K$ pierce the disk
bounded by $\eta$. Let $K_0=K(\eta,J)$ be the knot obtained
by replacing the $m$ trivial strands of $K$ by $m$ strands
``tied into the knot $J$''. More precisely, replace them with
$m$ untwisted parallels of a knotted arc with oriented knot
type $J$ as in Figure~\ref{infection}. We call $K_0$ the {\it
result of infecting $K$ by $J$ along $\eta$}.

\begin{figure}[ht!]
\cl{\relabelbox\small
\epsfxsize 2in \epsfbox{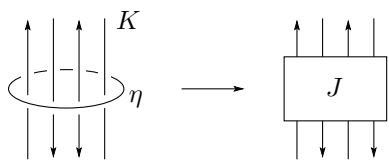}
\adjustrelabel <0pt, 0pt> {K}{$K$}
\adjustrelabel <0pt, 4pt> {e}{$\eta$}
\adjustrelabel <-2pt, 2pt> {J}{$J$}
\endrelabelbox}
\caption{Infecting $K$ by $J$ along $\eta$} \label{infection}
\end{figure}

The more general procedure of replacing the $m$ strands by a more
complicated string link will be discussed briefly in
Section~\ref{bordism}. Note that this is just a satellite
construction and as such is not new. The emphasis here is on
choosing the loop or loops $\eta$ to be very subtle with respect
to some measure. Note that this construction is, in a sense,
orthogonal to techniques used by Casson--Gordon, Litherland,
Gilmer, T.Stanford, and K.Habiro wherein the loop $\eta$ is
arbitrary but the analogue of the infection parameter $J$ is
increasingly subtle (for example, in Stanford's case, $J$ must lie
in the $n^{\supth}$ term of the lower central series of the pure
braid group; and, in the {\em claspers} that Habiro associated to
Vassiliev theory, the analogue of $\eta$ is a meridian of $K$
\cite{Hb}). However, infection can certainly be viewed as the
result of modifying $K$ by a certain {\it clasper} (depending on
$J$) all of whose {\it leaves} are parallels of $\eta$ (see
\cite{CT}\cite{GL}\cite{GR}). Moreover all of these procedures are
special cases of the classical technique, used by J. Levine and
others, of modifying a knot by Dehn surgeries that leave the
ambient manifold unchanged.

We now give an alternate description of genetic infection that is
better suited to analysis by Mayer--Vietoris and Seifert--Van Kampen
techniques. Beginning with the exterior of $K$, $E(K)$, delete the
interior of a tubular neighborhood $N$ of $\eta$ and replace it
with the exterior of $J$, $E(J)$, identifying the meridian
$\mu_\eta$ of $\eta$ with the longitude $\ell_J$ of $J$, and the
longitude $\ell_\eta$ of $\eta$ with the meridian $\mu_J$ of $J$.
It is well-known and is a good exercise for the reader to show
that the resulting space is $E(K_0)$ as described above. Note that
this replaces the exterior of a unknot with the exterior of the
knot $J$ in a fashion that preserves homology. Since there is a
degree one map (rel boundary) from $E(J)$ to $E$(unknot), there is
a degree one map (rel boundary) $f$ from $E(K_0)$ to $E(K)$ which
is the identity outside $E(J)$.

\begin{thm}\label{subtle} If $\eta\in\gs n$ then the map $f$
(above) induces an isomorphism
$f\co\pi_1(E(K_0))/\pi_1(E(K_0))^{(n+1)}\to\pi_1(E(K))/\pi_1(E(K))^{(n+1)}$
and hence induces isomorphisms between the $i^\supth$
(integral {\it and} localized) modules of $K_0$ and $K$ for
$0\le i<n$.
\end{thm}

\begin{proof} Let $E(\eta)$ denote $E(K)$ with the interior
of an open tubular neighborhood of $\eta$ deleted. Then, by the
Seifert--VanKampen theorem, $G=\pi_1(E(K))\cong\linebreak
\langle\pi_1(E(\eta)), t\mid\mu_\eta{=}1,\ell_\eta{=}t\rangle$.
Similarly,
$G_0=\pi_1(E(K_0))\cong\langle\pi_1(E(\eta)),\pi_1(E(J))\linebreak
|\ell_\eta=\mu_J,\mu_\eta=\ell_J\rangle$ where this denotes the
obvious ``free product with amalgamation''. The map $f$ induces
the identity on $\pi_1(E(\eta))$ and is the Hurewicz map on
$\pi_1(E(J))\to\BZ=\<t\>$ which sends $\ell_J\to1$ and $\mu_J\to
t$. Hence the kernel of $f\co G_0\to G$ is precisely the normal
closure in $G_0$ of $[P,P]$ where $P=\pi_1(E(J))$. Thus it
suffices to show that $P\subset\gs n_0$. Since $P$ is normally
generated by $\mu_J$, it suffices to show by induction that
$\mu_J\in\gs n_0$. This is clearly true for $n=0$. Suppose
$\mu_J\in\gs k_0$ $k<n$. Then $P\subset\gs k_0$ so
$\mu_\eta=\ell_J\subset[P,P]\subset G^{(k+1)}_0$. By hypothesis
$\eta\in\gs n$. Therefore $\eta$ bounds in $E(K)$, a map of a
symmetric $n$--stage grope \cite{CTe}. Thus $\ell_\eta$ bounds
such a grope in $E(K)$ and we may assume that the grope stages
meet $\eta$ transversely. Hence $\ell_\eta$ bounds a {\it
punctured} $n$--stage grope in $E(\eta)$ and the boundaries of
these punctures are copies of $\mu_\eta$. Therefore, in $G_0$,
$\ell_\eta=\prod^m_{i=1}\xi_i\mu^{n_i}_\eta\xi^{-1}_i\prod^r_{j=1}[a_j,b_j]$
where each $a_j$ and $b_j$ bound maps of punctured $(n-1)$--stage
gropes in $E(\eta)$. We claim $\ell_\eta\in G^{(k+1)}_0$. It
suffices to show the $a_j$ and $b_j$ lie in $\gs k_0$. But each of
these, modulo conjugates of $\mu^{\pm1}_\eta$, is given by a
similar expression as $\ell_\eta$ above. Continuing in this
fashion, we see that $\ell_\eta\in\gs n_0$ modulo the punctures
$\mu_\eta\in G^{(k+1)}_0$. Since $n\ge k+1$, $\ell_\eta\in
G^{(k+1)}_0$ and hence $\mu_J\in G^{(k+1)}_0$, completing our
induction.
\end{proof}

\begin{thm}\label{thm:infection} Let $K_0=K(\eta,J)$ be the result
of genetic infection of $K$ by $J$ along $\eta\in\gs n$ (as
described above). Then the $n^\supth$ (integral) Alexander
module of $K_0$, $\SA^\bbz_n(K_0)$, is isomorphic to
$\SA^\bbz_n(K)\op(\SA^\bbz_0(J)\ox_{\bbz[t,t^{-1}]}\bbz[G/\gns+])$
where $\bbz[G/\gns+]$ is a left $\bbz[t,t^{-1}]$ module via
the homomorphism $\<t\>=\bbz\to G/\gns+$ sending $t\to\eta$.
Thus, if $n\ge1$, $\SA_i(K_0)\cong\SA_i(K)$ for all $i\le n$.
\end{thm}

\begin{proof} Note that since $\eta\in\gs n$,
$\SA^\bbz_n(K_0)$ and $\SA^\bbz_n(K)$ are modules over
isomorphic rings since $G/\gns+\cong G_0/(G_0)^{(n+1)}$ by
the previous theorem. Therefore we can take the point of view
that the map $E(K_0)\to E(K)$ induces on both spaces a local
coefficient system with $G/\gns+$ coefficients.

\begin{lem}\label{boundary} The inclusion $i\co\p E(J)\to E(J)$
induces an isomorphism on \linebreak
$H_0(\un{\ \ };\bbz[G/\gns+])$ and
induces either the $0$ map or an epimorphism on \linebreak
$H_1(\un{\ \ };\bbz[G/\gns+]$ according
as $\eta\notin\gns+$ or $\eta\in\gns+$ respectively, whose
kernel is generated by $\<\ell_J\>$.
\end{lem}

\begin{proof}[Proof of Lemma~\ref{boundary}] The Lemma refers
to the coefficient system on $E(J)$ induced by $E(J)\sbq
E(K_0)\to E(K)$. Note that the kernel of the map
$\pi_1(E(J))\to G$ contains $[\pi_1(E(J)),\pi_1(E(J))]$
and thus its image in $G/\gns+$ is cyclic, generated by the
image of $\mu_J=\eta$. Since $G/\gns+$ is torsion-free (see
Example 2.4), this image is either zero or $\bbz$ according as
$\eta\in\gns+$ or not. This also shows that the image of
$\pi_1(\p E(J))$ in $G/\gns+$ is the same as the image of
$\pi_1(E(J))$. The first claim of the Lemma now follows
immediately from the proof of Proposition~\ref{Ho}.
Alternatively, since $H_0(\un{\ \ };\bbz[G/\gns+])$ is free on
the path components of the induced cover, and since the
cardinality of such is the index of the image of $\pi_1$ in
$G/\gns+$, $i$ induces an isomorphism on $H_0(\un{\ \
};\bbz[G/\gns+])$. If $\eta\in\gns+$ then the induced local
coefficient systems on $\p E(J)$ and $E(J)$ are trivial,
i.e.\ untwisted and thus $i$ induces an epimorphism on
$H_1(\un{\ \ };\bbz[G/\gns+])$ whose kernel is $\<\ell_J\>$
because it does so with ordinary $\bbz$ coefficients. If
$\eta\notin\gns+$ then the induced cover of $\p E(J)$ is a
disjoint union of copies of the $\bbz$--cover which ``unwinds''
$\mu_J$, i.e.\ the ordinary infinite cyclic cover. Thus $H_1$
of this cover is generated by a lift of $\ell_J$. But $\ell_J$
bounds a surface in $E(J)$ and this surface lifts to the
induced cover since every loop on a Seifert surface lies in
$[\pi_1(E(J)),\pi_1(E(J))]$. Therefore $i$ induces the zero
map on $H_1$ in this case. This concludes the proof of the
Lemma.
\end{proof}

We return to the proof of Theorem~\ref{thm:infection}. Consider
the Mayer--Vietoris sequence with $\bbz[G/\gns+]$ coefficients for
$E(K_0)$ viewed as $E(J)\cup E(\eta)$ with intersection $\p E(J)$.
By Lemma~\ref{boundary} this simplifies to
$$
H_1(\p E(J))\overset{(\psi_1,\psi_2)}{\lra} H_1(E(J))\op
H_1(E(\eta))\lra H_1(E(K_0))\lra0.
$$
Note first that $E(K)$ is obtained from $E(\eta)$ by adding a
solid torus, i.e.\ a 2-cell and then a 3-cell, so that it is
clear that $H_1(E(K))$ is the quotient of $H_1(E(\eta))$ by
the submodule generated by $\mu_\eta$ (or $\ell_J$). If
$\eta\notin\gns+$ then $\psi_1$ is zero by
Lemma~\ref{boundary} so $H_1(E(K_0))\cong
H_1(E(J))\op(H_1(E(\eta))/\<\psi_2\>)$. But in the proof of
Lemma~\ref{boundary} we saw that $H_1(\p E(J))$ was generated
by $\ell_J$ and so the image of $\psi_2$ is generated by
$\ell_J$. Hence $H_1(E(K_0))\cong H_1(E(J))\op(H_1(E(K))$.
This concludes the proof of the theorem in the case
$\eta\notin\gns+$ once we identify $H_1(E(J))$ as
$\SA^\bbz_0(J)\ox_{\bbz[t,t^{-1}]}\bbz[G/\gns+]$. But since
the map from $\pi_1(E(J))$ to its image in $G/\gns+$ has
already been observed to be the abelianization, this is clear.

In case $\eta\in\gns+$, $\psi_1$ is an epimorphism whose
kernel is generated by $\ell_J$ and so $H_1(E(K_0))\cong
H_1(E(\eta))/\<\ell_J\>\cong H_1(E(K))$. On the other hand,
in this case $\SA^\bbz_0(J)\ox_{\bbz[t,t^{-1}]}\bbz[G/\gns+]$
factors through the augmentation of $\SA^\bbz_0(J)$, which is
zero since the classical Alexander polynomial of a knot
augments to $1$.

If $n\ge1$, $\eta\in\wt G$. Let $\Delta(t)$ be the classical
Alexander polynomial of $J$. Then $\Delta(\eta)\in\bbz\wt
G-\{0\}$. Recall that $\bbz\wt G-\{0\}$ is a right divisor
set of regular elements of $\bbz[G/\gns+]$ by
\cite[p.~609]{P}. Thus for any $r\in\bbz[G/\gns+]$, there
exist $r_1\in\bbz[G/\gns+]$ and $t_1\in\bbz\wt G-\{0\}$ such
that $\Delta(\eta)r_1=rt_1$ \cite[p.~427]{P}. Hence any
element $x\ox r\in\SA^\bbz_0(J)\ox\bbz[G/\gns+]$ is
annihilated by $t_1$, showing that this is a $\bbz\wt
G$--torsion module. Hence $\SA_n(K_0)\cong\SA_n(K)$.
\end{proof}

\section{Applications to detecting fibered and alternating
knots and symplectic structures on
$4$--manifolds}\label{fiber}

In this section we show that the higher-order Alexander modules of
fibered knots and alternating knots have special properties.
Therefore noncommutative knot theory gives algebraic invariants
which can be used to tell when a knot is not fibered or not
alternating, even in situations where the Alexander module yields
inconclusive evidence. In the case of fibered knots, examples of
this type were obtained independently by J.C. Cha using the
twisted Alexander invariant \cite{Ch}. Remarkably, for
$4$--manifolds of the form $M_K\x S^1$ ($M_K$ is the $0$--surgery on
$K$), our invariants also obstruct the existence of a symplectic
structure (using work of P. Kronheimer \cite{Kr}). We also
establish that $\d_i-\d_j$ are not Vassiliev invariants of finite
type.

\begin{prop}\label{fibered} If $K$ is a non-trivial fibered
or alternating knot then
$\d_0=\d_1+1=\dots=\d_n+1=2\genus(K)$.
\end{prop}

\begin{proof} It is well known that for a fibered or
alternating knot, $\d_0=2\genus(K)$.  The result now follows
immediately from Theorem~\ref{non-decreasing}.
\end{proof}

\begin{cor}\label{fibered example} For any non-trivial fibered
or alternating knot $K$, and any positive integer $n$, there
exists a hyperbolic knot $K_*$ such that
\begin{description}
\item[a)] $\SA_k(K_*)\cong\SA_k(K)$ for all $k<n$
\item[b)] $\d_k(K_*)=\d_k(K)$ for all $k<n$
\item[c)] $\d_n(K^*)>\d_n(K)$
\item[d)] $K_*$ is neither fibered nor alternating.
\end{description}
\end{cor}

\begin{proof} Apply Theorem~\ref{increasing degree} and
Corollary~\ref{hyperbolic} to produce $K_*$. Suppose $n\ge2$.
Since $K$ is fibered or alternating, $\d_n(K)=\d_{n-1}(K)$ by
Proposition~\ref{fibered}. It follows that
$\d_n(K_*)>\d_{n-1}(K_*)$ so $K_*$ is not fibered. A similar
argument works for $n=1$.
\end{proof}

There are more subtle obstructions to fibering that cannot be
detected by the localized modules, but can be detected by the
integral modules.

\begin{prop}\label{torsionfree} If $K$ is a fibered knot then
the following equivalent conditions hold:
\begin{enumerate}
\item[{\rm1)}] $\SA^\bbz_n(K)\lra\SA_n(K)$ is injective
\item[{\rm2)}] $\SA^\bbz_n(K)$ is torsion-free as a $\bbz\wt
G$--module (recall that $\wt G$ is $G^{(1)}/\gns+$).
\end{enumerate}
\end{prop}

\begin{proof} Recall
$\SA^\bbz_n(K)=\gns+/G^{(n+2)}=F^{(n)}/F^{(n+1)}$ where
$G^{(1)}=F$ is free since $K$ is a fibered knot. Since $\wt
G=G^{(1)}/\gns+=F/F^{(n)}$, $\SA^\bbz_n$ as a $\bbz\wt
G$--module is merely $F^{(n)}/F^{(n+1)}$ as a
$\bbz[F/F^{(n)}]$--module (i.e.\ $H_1(F;\bbz[F/F^{(n)}])$).
Since $F$ is the fundamental group of a $1$--complex, this is
a submodule of a free module and hence is torsion-free.
\end{proof}

\begin{thm}\label{subtlefibered} For any non-trivial fibered
knot $K$ and any positive integer $n$ there exists a family
of hyperbolic knots $K_*=K_*(J,n)$, parametrized by an
auxiliary knot $J$, such that
\begin{enumerate}
\item[{\rm1)}] $G/\gns+\cong G_*/\gns+_*$ meaning that all knots
in the family share (with $K$) the same $\SA^\bbz_i$ for $0\le
i\le n-1$; \item[{\rm2)}] $\SA_n(K)\cong\SA_n(K_*)$ \item[{\rm3)}]
$\d_0=\d_1+1=\dots=\d_n+1$ for each $K_*$ and $K$ \item[{\rm4)}]
if $P_*$ is the commutator subgroup of $G_*$ then
$P_*/(P_*)_j\cong F/F_j$ for each term of the lower central series
($F$ is free of rank equal to $2\genus(K)$). \item[{\rm5)}] If $J$
has non-trivial classical Alexander polynomial then $K_*$ is not
fibered and hence is distinct from $K$. \item[{\rm6)}] If $J$ has
non-trivial classical Alexander polynomial then
$G/G^{(n+2)}\not\cong$ $ G_*/G^{(n+2)}_*$ and
$\SA^\bbz_n(K_*)\not\cong\SA^\bbz_n(K)$.
 \item[{\rm7)}] $K_*(J,n)$ and $K_*(J',n)$ are distinct if the integrals of the classical
Levine signature functions of $J$ and $J'$ are distinct.
\end{enumerate}
\end{thm}

\begin{proof} Since $K$ is a non-trivial fibered knot, it
does not have Alexander polynomial $1$. By
Corollary~\ref{non-triviality} (or more simply since $\gs n$ is
free if $n\ge1$), for any $n$, we can choose a class $\eta\in\gs
n-\gns+$ which can be represented by a loop in the complement of a
fiber surface for $K$ and which is also unknotted in $S^3$.
Construct $K_0=K(\eta,J)$ by genetic infection as in
Section~\ref{genetic}. It follows that $\genus(K_0)=\genus(K)$ so
if we drop the claim of hyperbolicity we can retain this. By
Corollary~\ref{hyperbolic} there is a hyperbolic knot, $K_*$,
whose fundamental group differs from that of $K_0$ by a perfect
group. Thus $K_0$ and $K_*$ have isomorphic $\SA^\bbz_i$ for any
$i$ and have isomorphic groups modulo any term of the derived
series (see Proposition~\ref{derived}). Thus, by
Theorem~\ref{subtle}, part 1) of Theorem~\ref{subtlefibered}
follows. Part 2) follows from Theorem~\ref{thm:infection}. Part 3)
holds for $K$ by Proposition~\ref{fibered} and hence for $K_*$ by
the second part of Theorem~\ref{thm:infection}. Part~4) is true
for {\it any} knot for which $\SA^\bbz_0\cong\bbz^{2\genus(K)}$
since one can then define a homomorphism from the free group of
rank $2\genus(K)$ to the commutator subgroup which induces an
isomorphism on $H_1$ and an epimorphism on $H_2$. Stallings'
theorem \cite{St} then guarantees an isomorphism modulo any term
of the lower central series.

By Theorem~\ref{thm:infection}, $\SA^\bbz_n(K_*)$ is
$\SA^\bbz_n(K)$ direct sum
$\SA^\bbz_0(J)\ox_{\bbz[t,t^{-1}]}\bbz[G/\gns+]$. But
$\SA_n(K_*)\cong\SA_n(K)$. By Proposition~\ref{torsionfree}, if
$K_*$ were fibered then this second direct summand would be zero.
But even after tensoring with $\bbq[G/\gns+]$ this module is not
zero because it is cyclic of order $\Delta(\eta)$ where
$\Delta(t)$ is the classical Alexander polynomial of $J$. Thus the
module is zero if and only if $\Delta(\eta)$ is a unit in
$\bbq[G/\gns+]$. Since $G/\gns+$ is PTFA it is right orderable by
\cite[p.~587]{P} hence has only trivial units by
\cite[p.~588,590]{P}. Since $\Delta(\eta)$ is an integral
polynomial in $\eta$, this can only happen if $\Delta(t)$ has
degree zero which was excluding by hypothesis. Thus Part 5 is
established. Part 6 follows from the discussion above. The proof
of Part 7) must be postponed to Theorem~\ref{subtle signature} of
Section~\ref{signatures}.
\end{proof}

Some of these new obstructions to fibering can be used to
show that certain 4-manifolds of the form $S^1\x M_K$ admit
no symplectic structure. If $K$ is a fibered knot then $M_K$
also fibers over the circle and it is known that $S^1\x M_K$
is then symplectic. C. Taubes conjectured the converse. The
Seiberg-Witten invariants provide evidence for this
conjecture. If $S^1\x M_K$ admits a symplectic structure then
the Alexander polynomial of $K$ must be monic. This is
precisely the fibering obstruction on the classical Alexander
polynomial of $K$. Peter Kronheimer provided more evidence
for the conjecture by proving that if $S^1\x M_K$ admits a
symplectic structure then $\d_0=2\genus(K)$
\cite{Kr2} \cite{Kr}. As a consequence of his work, we see
that the $\d_i$ constitute algebraic invariants which can
obstruct a symplectic structure on $S^1\x M_K$ even when the
Seiberg Witten invariants give inconclusive information.

\begin{thm}\label{symplectic} Suppose $K$ is a non-trivial
knot. If $S^1\x M_K$ admits a symplectic structure then the
invariants $\d_i(K)-\d_0(K)+1$ are zero for all $i>0$.
\end{thm}

\begin{proof} By Kronheimer's theorem, $\d_0(K)=2\genus(K)$.
The result then follows from Theorem~\ref{non-decreasing}.
\end{proof}

\begin{cor}\label{symplectic examples} If $K_*$ is any one
of the examples of Corollary~\ref{fibered example}, then
$S^1\x M_K$ admits no symplectic structure although the
Alexander polynomial of $K$ is monic.
\end{cor}

Now consider $\d_n$ as a rational valued invariant on knot
types.

\begin{prop}\label{finite type} None of the invariants
$\d_i-\d_j$ ($i\neq j$) or $\d_i-2\genus(K)$ is determined by
any finite number of finite type (Vassilliev) invariants.
\end{prop}

\begin{proof} Let $\phi$ be one of the mentioned invariants.
Suppose $\phi$ were determined by the finite type invariants
$v_1,\dots,v_m$. We have shown in Theorem~\ref{increasing
degree} that $\phi$ is not constant. But on fibered knots
$\phi$ is constant, say $C$, by Proposition~\ref{fibered}.
Again using Theorem~\ref{increasing degree} choose a knot
$K_*$ such that $\phi(K_*)=C'\neq C$. By a result of A.
Stoimenow \cite{Sti}, there exists a fibered knot $K$ such
that $v_i(K)=v_i(K_*)$ for $1\le i\le m$. Thus
$\phi(K)=\phi(K_*)$ a contradiction.
\end{proof}

\section{Bordism invariants generalizing the Arf
invariant}\label{bordism} In this section we define
higher-order bordism invariants for knots which (in a certain
sense) generalize the Arf invariant. The reader is warned that
these are not the same as the generalizations of the Arf
invariant defined in Section 4 of \cite{COT1}. The invariants
about to be defined are almost certainly {\it not} concordance
invariants. If $K$ is a knot, $G$ its group, let $M_K$ be the
result of $0$--framed surgery on $K$ and $P=\pi_1(M_K)$. Recall
that the Arf invariant of $K$ may be defined as the class in
$\Om^\spin_3(S^1)\cong\bbz_2$ represented by $M_K$ with the
map to $S^1$ induced by the abelianization homomorphism
$P\lra P/P^{(1)}\cong G/\gs1\cong\bbz$. Equivalently one
could consider spin bordism (rel boundary) over $S^1$ of
$3$--manifolds with a toral boundary component, in which case
the Arf invariant of $K$ is zero if and only if $\sk$ is
bordant to the exterior of the unknot. Note that
$S^1=K(P/P^{(1)},1)=K(G/\gs1,1)$.

More generally,

\begin{defn}\label{reduced bordism} The {\it $n^{\rm th}$
(reduced) bordism invariant of $K$}, $\b_n(K)$, is the class
in $\Om^\spin_3(K(P/P^{(n+1)},1))/\aut(P/P^{(n+1)})$
represented by $M_K\overset{f_n}{\lra}\linebreak
K(P/P^{(n+1)},1)$ where
$f_n$ is induced by the quotient map $f_n\co P\lra P/P^{(n+1)}$.
\end{defn}

Obviously then $\b_0(K)$ is the Arf invariant of $K$, and
equivalent knots (here we need an {\it orientation-preserving}
homeomorphism) have identical bordism invariants. Note also that
$\b_n(-K)=-\b_n(K)$, so that a $\pm$--amphichiral knot satisfies
$2\b_n=0$. We also have the following purer but uglier version.
The purest (and ugliest) version would fix the peripheral
structure in $G/\gns+$ and effectively only consider pairs of
knots with isomorphic $G/\gns+$ preserving peripheral structure.

\begin{defn}\label{unreduced bordism} The {\it $n^{\rm th}$
unreduced bordism invariant of $K$}, $\tl\b_n(K)$, is the
equivalence class of ($\sk,\p(\sk)$, $f_n\co\sk\lra
K(G/\gns+,1))$ in the set of spin bordism classes rel
boundary of spin $3$--manifolds and maps to $K(G/\gns+,1)$
modulo the action of $\aut(G/\gns+,1)$.
\end{defn}

\begin{con} For each $n\ge0$ there exist knots $K$, $K_*$
such that $\SA_i(K)\cong\SA_i(K_*)$, $0\le i<n$ but
$\b_n(K)\neq\b_n(K_*)$.
\end{con}

We give a construction which should produce $K_*$ for any $n$
but are only able to verify this in the first case $n=0$.

\begin{thm}\label{bordism theorem} There exist knots $K$,
$K_*$ with identical classical Alexander module and Blanchfield
form but which are distinguished by $\b_1$. Moreover we can choose
$K$ to be amphichiral and $K_*$ to be chiral.
\end{thm}

The knot $K_*$ in this case will be constructed from $K$ by
choosing 3 ``bands'' of a Seifert surface for $K$ and tying
them into the shape of a Borromean rings, with a restriction
on the 3 bands that they are ``essential'' (in a sense to be
made precise) in the Alexander module. This then becomes
very interesting in light of previous work of S. Naik and T.
Stanford who showed that any two knots with isomorphic
classical Alexander modules and isomorphic classical
Blanchfield forms are related by a sequence of such
replacements (without the restriction) \cite{NS}; and work of
S. Garoufalidis and L. Rozansky, who define a ``finite type''
isotopy invariant of knots that is affected by precisely this
same construction. Their invariant gives information even for
Alexander polynomial one knots \cite{GR}. The work of
Naik-Stanford can be interpreted as saying that the
construction by which we prove Theorem~\ref{bordism theorem}
is the {\it only} construction necessary (for $n=0$) to
achieve the full range of values of the triple ($\SA_0(K)$,
$\SB\ell_0(K)$, $\b_1(K)$).

\medskip
\sh{Infection by a string link}

\begin{figure}[ht!]
\cl{\relabelbox\small
\epsfxsize 2.3in \epsfbox{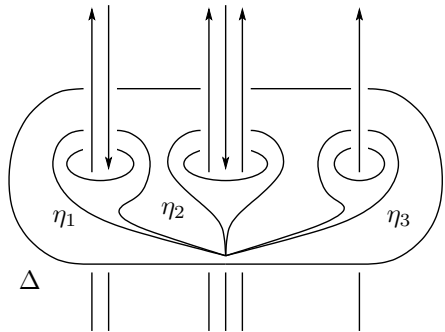}
\adjustrelabel <0pt, 0pt> {D}{$\Delta$}
\adjustrelabel <0pt, 0pt> {e}{$\eta_1$}
\adjustrelabel <0pt, 0pt> {f}{$\eta_2$}
\adjustrelabel <0pt, 0pt> {g}{$\eta_3$}
\endrelabelbox}
\caption{The data required to infect $K$ by a string link
along $(\eta_1,\eta_2,\eta_3)$} \label{Borromean}
\end{figure}

We discuss an instance of genetic infection of a fixed knot $K$
using an auxiliary {\it string link}. This was perhaps first
discussed in \cite{CO} in the case that the auxiliary link is a
boundary link. For our examples it suffices to consider the case
where the auxiliary link is the Borromean Rings. This type of
Borromean modification has been considered by many other authors
including Matveev, Habiro, Goussarov and those mentioned above.
Let $\Delta$ be a 2-disk with 3 disjoint open subdisks $\Delta_1$,
$\Delta_2$, $\Delta_3$, deleted. Consider an embedding of $\Delta$
in $\sk$ which extends to an embedding $\Delta^+$ of $D^2$ into
$S^3$. An example is shown in Figure~\ref{Borromean}. The trivial
braid $K\cap(\Delta^+\x[0,1])\hra\Delta^+\x[0,1]$ is obtained from
the trivial 3-string braid by forming $\{m_1,m_2,m_3\}$ parallel
strands (and perhaps altering some orientations) where $m_i$ is
the number of components of $K\cap\Delta_i$. From the Borromean
rings (written as a 3-string braid) form the $\{m_1,m_2,m_3\}$
cable of the Borromean rings and alter orientations consistent
with the above. Then replace the trivial braid with this cable of
the Borromean Rings. We denote the modified knot by $K_*=K(\eta)$
where $\eta$ denotes the triple
$(\eta_1,\eta_2,\eta_3)=(\p\Delta_1,\p\Delta_2,\p\Delta_3)$ of
conjugacy classes of elements of $G=\pi_1(E(K))$ as shown in
Figure~\ref{Borromean}.

Once again this is the same as replacing the solid handlebody
$\Delta\x[0,1]$ with the exterior of a 3-string braid that
represents the Borromean Rings.

The Seifert Van-Kampen proof of the following Lemma is left
to the reader, it being entirely analogous to that of
Theorem~\ref{subtle}.

\begin{lem}\label{Borromean subtle} If $\eta_i\in\gs n$,
$i=1$, $2$, $3$ then $\pi_1(E(K_*))/(\pi_1(E(K_*)))^{(n+1)}$
is isomorphic to $\pi_1(E(K))/(\pi_1(E(K)))^{(n+1)}$
preserving peripheral structure. In particular
$\SA^\bbz_i(K_*)\cong\SA^\bbz_i(K)$ for $0\le i<n$.
\end{lem}

\begin{figure}[ht!]
\cl{\relabelbox\small
\epsfxsize 2.3in \epsfbox{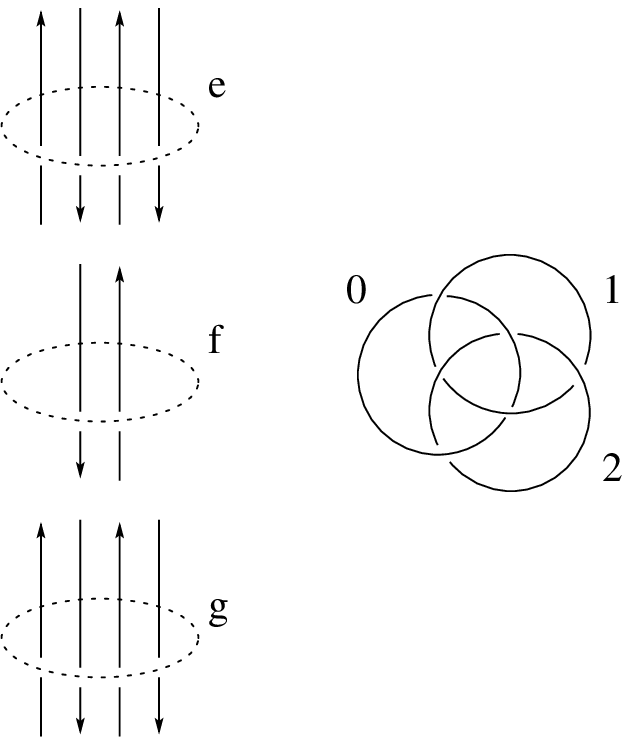}
\adjustrelabel <0pt, 0pt> {0}{$0$}
\adjustrelabel <0pt, 0pt> {1}{$0$}
\adjustrelabel <0pt, 0pt> {2}{$0$}
\adjustrelabel <0pt, 0pt> {e}{$\eta_1$}
\adjustrelabel <0pt, 0pt> {f}{$\eta_2$}
\adjustrelabel <0pt, 0pt> {g}{$\eta_3$}
\endrelabelbox}
\nocolon\caption{} \label{cobordism}
\end{figure}

\begin{figure}[ht!]
\cl{\relabelbox\small
\epsfxsize 2.3in \epsfbox{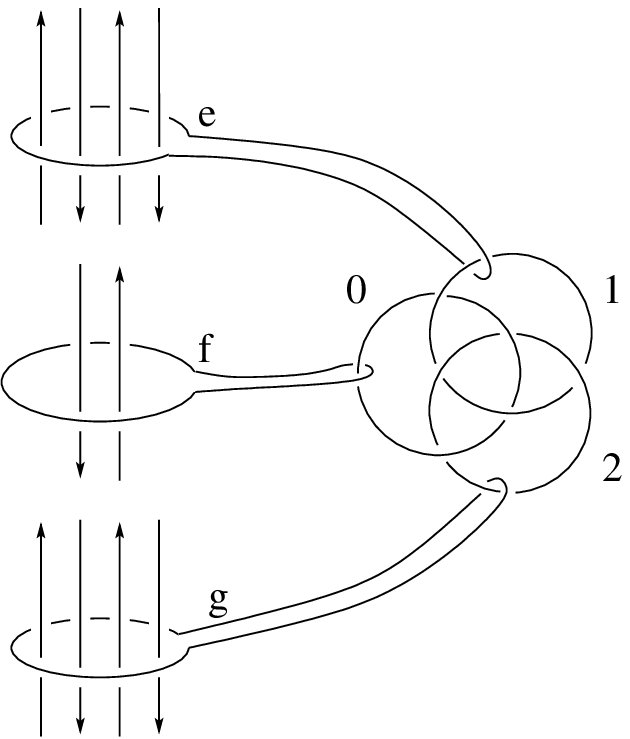}
\adjustrelabel <0pt, 0pt> {0}{$0$}
\adjustrelabel <0pt, 0pt> {1}{$0$}
\adjustrelabel <0pt, 0pt> {2}{$0$}
\adjustrelabel <0pt, 0pt> {e}{$0$}
\adjustrelabel <0pt, 0pt> {f}{$0$}
\adjustrelabel <0pt, 0pt> {g}{$0$}
\endrelabelbox}
\nocolon\caption{} \label{cobordism2}
\end{figure}

Now consider the difference $\b_n(K_*)-\b_n(K)$ projected onto
$\Om_3(P/P^{(n+1)})\cong H_3(P/P^{(n+1)})$, forgetting the
spin structure. We claim that this element is equal to the
image of a generator of $H_3(\sss)$ under the map induced by
$\bbz^3\overset{(\eta_1,\eta_2,\eta_3)}{\lra}P^{(n)}/P^{(n+1)}\lra
P/P^{(n+1)}$. To establish this, we describe a cobordism, over
$P/P^{(n+1)}$, from
($M_K\amalg\sss,f_n\amalg(\eta_1,\eta_2,\eta_3)$) to
$(M_{K_*},(f_*)_n)$. First add a 1-handle to $\p_+$ of
$(M_K\amalg(\sss))\x[0,1])$. A framed link picture of the new
$\p_+$ is shown in Figure~\ref{cobordism}. Since the
meridians of the components of the pictured Borromean rings
map to $(\eta_1,\eta_2,\eta_3)$ in $P/P^{(n+1)}$ we can add
three 2-handles $\{h_1,h_2,h_3\}$ as shown in
Figure~\ref{cobordism2} and still have a cobordism over
$P/P^{(n+1)}$. But now the knot in Figure~\ref{cobordism2} is
well known to be equivalent to $K(\eta)$ by first sliding the
strands of $K$ which link the attaching circles of $h_i$ over
the corresponding component of the Borromean Rings until the
attaching circles of the $h_i$ bound disks intersecting only
the Borromean rings and then sliding the strands of $K$ over
the $h_i$ as needed until completely free.

\begin{proof}[Proof of Theorem~\ref{bordism theorem}]
Let $K$ be a knot with classical Alexander module cyclic
 of order $p(t)p(t^{-1})$ where $p(t)=t^3+t-1$. Since
$p(t)$ is irreducible and coprime to $p(t^{-1})$ there is a unique
direct summand $B$ of $\SA_0(K)$ isomorphic to
$\bbz[t,t^{-1}]/\<p(t)\>$. Since $B$ is a free abelian group of
rank $3$, $1\land t\land t^2$ is a basis of $H_3(B)$ and also
represents an element $\a$ of $H_3(G/\gs2)\cong H_3(P/P^{(2)})$
under the inclusions. Choose a trivial link
$\{\eta_1,\eta_2,\eta_3\}$ in $S^3$, avoiding $K$, representing
$\{1,t,t^2\}$ in $G/\gs2$ and perform a Borromean modification to
$K$ along $\{\eta_1,\eta_2,\eta_3\}$ as above to arrive at a knot
$K_*$ that has isomorphic $\SA_0$. Using the cobordism above we
see that $\b_1(K_*)-\b_1(K)=\a$. However we must take into account
the ambiguity in the definition of $\b_1(K)$. Suppose $f$ is an
automorphism of the group $G/\gs2$. Assuming
$\b_1(K_*)-f_*\b_1(K)=0$ for some $f$, we shall derive a
contradiction. Here we are viewing both $\b_1(K_*)$ and $\b_1(K)$
as elements of $H_3(G/\gs2)$. Let $r$ be the canonical retract
from $\bbz[t,t^{-1}]/\<p(t)p(t^{-1})\>$ to $B$, inducing a map $r$
from $G/\gs2$ (which is $\SA_0\rtimes\bbz$ to $\pi=B\rtimes\bbz$).
The automorphism $f$ induces an automorphism $g$ of $\pi$ such
that $g\circ r=r\circ f$. Combining the two equations above we see
that $f_*\b_1(K)-\b_1(K)=\a$ and hence
$r_*(\a)=r_*f_*\b_1(K)-r_*\b_1(K)=(g_*-\id)(r_*\b_1(K))$. Consider
the Wang sequence
$$
H_3(B)\overset{t_*-\id}{\lra}H_3(B)\overset{j_*}{\lra}H_3(\pi)
\overset{\p}{\lra}H_2(B)\overset{t_*-\id}{\lra}H_2(B).
$$
Since $H_2(B)$ is free abelian on $\{1\land t,1\land t^2,t\land
t^2\}$ one can easily calculate that $(t_*-\id)$ is injective on
$H_2(B)$. Hence $\p$ is the zero map and $r_*(\b_1(K))=j_*(\b)$
for some $\b\in H_3(B)$. Recall that, by definition,
$r_*(\a)=j_*(1\land t\land t^2)$. It follows that $(1\land t\land
t^2)-(g_*-\id)(\b)$ lies in the kernel of $j_*$ and hence in the
image of $(t_*-\id)$. But $H_3(B)$ is $\bbz$ generated by $1\land
t\land t^2$ so it is easy to calculate that $t_*-\id$ is zero on
$H_3(B)$. Moreover since $g_*$ is an automorphism of an infinite
cyclic group, it equals $\pm\id$. Hence $1\land t\land t^2=0$ or
$1\land t\land t^2=-2\b$, both contradictions. Therefore
$\b_1(K_*)$ and $\b_1(K)$ are distinct.
Alternatively, we could
choose the amphichiral knot $K\#-K$ and form $K_*$ by infecting
``the $K$ part'' as above. Then $K_*$ is not amphichiral.
\end{proof}

\section{Von Neumann higher-order signatures of
knots}\label{signatures}

One can define higher-order signatures $\rho_n$, $n\ge0$, for
knots using the Von Neumann $\rho$--invariant of J. Cheeger and M.
Gromov. In this section these are defined and used to distinguish
among knots which have isomorphic localized Alexander modules.
These can also be used to detect chirality of knots. Similar
signatures were crucial in the work of Cochran-Orr-Teichner
\cite{COT1} \cite{COT2} \cite{CT}.

If $K$ is a knot and $G$ its group, let $M_K$ denote the result of
zero framed surgery on $K$ and let $P=\pi_1(M_K)$. To the
$P^{(n+1)}$ covering space of $M_K$, Cheeger and Gromov associate
a real-valued {\it Von Neumann $\rho$--invariant}, which we denote
$\rho_n(K)$ \cite{ChG}. If $-K$ denotes the mirror-image of $K$
then $M_{-K}\cong-M_K$ so $\rho_n(-K)=-\rho_n(K)$ \cite{ChG}. If
$K$ and $K'$ are equivalent knots, then $M_K$ and $M_{K'}$ are
(orientation-preserving) homeomorphic so $\rho_n(K)=\rho_n(K')$.
Hence if $K$ is plus or minus amphichiral then $\rho_n(K)=0$ for
each $n$. In general it is not known how to compute $\rho_n$.
However {\it relative signatures} $\rho_n(K_0)-\rho_n(K_1)$ are
often easy to compute. Suppose $K_0$, $K_1$ are knots such that
$P_0/P^{(n+1)}_0\cong P_1/P^{(n+1)}_1$ where $P_i=\pi_1(M_{K_i})$
as above. Moreover suppose $M_{K_0}$ and $M_{K_1}$ are bordant
over $P/P^{(n+1)}$, as in the previous section, that is there
exists a compact oriented 4-manifold ($W$, $\psi\co\pi_1(W)\lra
P/P^{(n+1)}$) whose boundary is ($M_{K_0}$,
$\phi_0\co\pi_1(M_{K_0})\lra P/P^{(n+1)}$)
$\amalg(-M_{K_1},\phi_1\co\pi_1(M_{K_1})\lra P/P^{(n+1)})$ where
$\phi_i$ is a composition of the projection $P_i\lra
P_i/P^{(n+1)}_i$ with an arbitrary identification of
$P_i/P^{(n+1)}_i$ with a standard copy called $P/P^{(n+1)}$. Then
the relative signature $\rho_n(K_0)-\rho_n(K_1)$ is equal to the
(reduced) $L^2$--signature $\s^{(2)}_n(W)-\s(W)$ associated to
$\psi$ (see \cite[Section 5]{COT1}. This is often calculable. For
example, if $n=0$ and $K$ is an Arf invariant zero knot, then
$P/P^{(1)}\cong\bbz$ and it is known that ($M_K$,
$\phi_0\co P\lra\bbz$) is null-bordant, i.e.\ that $M_K$ bounds
($V,\psi$), and that $\s^{(2)}_0(K)$ is the {\it integral} of the
{\it Levine signature function} of $K$ and $\s(V)$ is the ordinary
knot signature \cite{COT2}. In this sense the $\rho_n$ generalize
``ordinary'' Levine-Tristram signatures associated to the
$\bbz$--cover of $\sk$.

The technique of {\it genetic infection} may be used to
modify a given knot $K$ so subtly that the two have isomorphic
$i^{\supth}$ localized modules for $i\le n$ and have
isomorphic (integral) modules for $i<n$. The difference in
the (integral) modules at the $n^\supth$ stage can (in many
cases) be detected by the $n^\supth$ relative signature.

\begin{thm}\label{subtle signature} Let $K_*=K(\eta,J)$ be
the result of genetic infection of $K$ by $J$ along
$\eta\in\gs n$ ($G=\pi_1(E(K))$) as in Section~\ref{genetic}.
Then
\begin{enumerate}
\item[{\rm1)}] $\SA^\bbz_i(K_*)\cong\SA^\bbz_i(K)$ for $i<n$,
\item[{\rm2)}] If $n>0$, $\SA_i(K_*)\cong\SA_i(K)$ and
$\d_i(K_*)=\d_i(K)$ for $i\le n$,
\item[{\rm3)}] If $\arf J=0$ then $\b_i(K_*)=\b_i(K)$ and
$\tl\b_i(K_*)=\tl\b_i(K)$ for $i\le n$,
\item[{\rm4)}] $\rho_i(K_*)=\rho_i(K)$ for $i<n$,
\item[{\rm5)}] If $\eta\notin P^{(n+1)}$ then
$\rho_n(K_*)-\rho_n(K)$ is the integral of the normalized
Levine signature function of $J$. If this real number is
non-zero then $K_*$ is distinct from $K$ and distinct from
the mirror image of $K_*$.
\end{enumerate}
\end{thm}

\begin{proof} Since 1) and 2) were shown in
Theorem~\ref{subtle} and Theorem~\ref{thm:infection} we begin
with 3). Consider the map $f\co M_J\lra S^1$ induced by the
abelianization of $\pi_1(E(J))$. Since $\Om_3(S^1)\cong0$ and
$\Om^\spin_3(S^1)\cong\bbz_2$ as detected by the Arf invariant of
$J$, it can be shown that $M_J$ is the boundary of a $4$--manifold
$V$ with $\pi_1(V)\cong\bbz$ generated by the meridian of $J$ and
such that $V$ extends the usual spin structure on $M_J$ if $\arf
J=0$. We may also assume signature$(V)=0$ by connected summing
with $\pm\bbc P(2)$'s. The boundary of $V$ decomposes into
$E(J)\cup(S^1\x D^2)$. We form a cobordism $W$ from $M_K$ to
$M_{K_*}$ (or $E(K)$ to $E(K_*)$) as follows. Let $W$ be the
$4$--manifold obtained from $M_K\x[0,1]$ by identifying $S^1\x
D^2\hra\p V$ with the solid torus neighborhood of $\eta$ in
$M_K\x\{1\}$ in such a way that $\p^+W=M_{K_*}$. Since
$\pi_1(M_K)\lra\pi_1(W)$ is an isomorphism, $W$ is a cobordism
``over'' $\pi_1(M_K)/(\pi_1(M_K))^{(i+1)}$ for any $i$. By
Theorem~\ref{subtle}, this quotient is isomorphic to that of
$K_*$ if $i\le n$ (since the longitudes are preserved under
the map $f$ of Theorem~\ref{subtle}) and so $\b_i$ and
$\tl\b_i$ agree for $K_*$ and $K$ if $i\le n$. Then
$\rho_n(K_*)-\rho_n(K)$ is equal to the $L^2$--signature of
$W$ associated the homomorphism
$\pi_1(W)\cong P\lra P/P^{(n+1)}$. By additivity of signature
\cite[Lemma 5.9]{COT1}, this is equal to the $L^2$ signature
of $V$ associated to the map $\pi_1(V)\lra\pi_1(W)\lra
P/P^{(n+1)}$. Since $\pi_1(V)\cong\bbz$ generated by
$\eta(=\mu_J)$, if $\eta\notin P^{(n+1)}$, this map is
injective. It follows that the $L^2$ signature associated to
$P/P^{(n+1)}$ is equal to that associated to the map
$\pi_1(V)\lra\bbz$ (its image) \cite[Proposition 5.13]{COT1}.
But this is the integral of the classical Levine signature
function of $J$ over the circle as remarked above
\cite[Appendix]{COT2}. Note that $\rho_i(K_*)=\rho_i(K)$ for
$i<n$ because in this case the map $\pi_1(V)\lra P/P^{(i+1)}$
is zero and the $L^2$--signature of $V$ is equal to its usual
signature (which is zero).
\end{proof}

\begin{rems} If $n\ge1$, it is easy to get
$\genus(K_*)=\genus(K)$ by choosing $\eta$ in the complement
of a minimal genus Seifert surface for $K$. Then $K$ and
$K_*$ also have identical Seifert form. This shows that
$\rho_n$ is determined neither by the localized modules or
$\d_i$ for $i\le n$, nor by the bordism invariants $\b_i$ for
$i\le n$, nor by the genus. Note that the above proof also
establishes part 7) of Theorem~\ref{subtlefibered}
\end{rems}

\begin{ques} Is $\rho_n(K)$ determined by $\SA_n(K)$ and the
$n^{\supth}$ linking form $B\ell_n$ discussed in the next
section?
\end{ques}

\section{Higher order Blanchfield linking forms, duality, and
the behavior of the longitude}\label{non-singular}

We will now show that the Blanchfield linking form defined on
the classical Alexander module generalizes to linking forms
$B\ell_n$ on the localized higher-order Alexander modules
$\SA_n$. We see that if $n\neq1$, we can get a
{\it non-singular} linking form. If $n=1$ the form is
non-singular after killing the longitude. Hence the $\SA_n$
are {\it self-dual} if $n\neq1$. Recall that $(\SA,\lambda)$
is a {\it symmetric linking form\/} if $\SA$ is a torsion
$\SR$--module and
$$
\lambda\co\SA\lra\ov{\Hom_\SR(\SA,\SK/\SR)}\equiv\SA^\#
$$
is an $\SR$--module map such that
$\lambda(x)(y)=\overline{\lambda(y)(x)}$ (here $\SK$ is the
field of fractions of $\SR$ and $\Hom$, which is naturally a
left module, is made into a right $\SR$--module using the
involution of $\SR$). The linking form is {\it
non-singular\/} if $\lambda$ is an isomorphism.

\begin{thm}\label{higher Bforms} {\rm\cite{COT1}}\qua Suppose $M$
is a compact, oriented, connected $3$--manifold with
$\b_1(M)=1$ and $\phi\co\pi_1(M)\lra\G$ a non-trivial PTFA
coefficient system. Suppose $\SR$ is a ring such that
$\bbz\G\sbq\SR\sbq\SK$. Then there is a symmetric linking form
$$
B\l\co H_1(M;\SR)\lra H_1(M;\SR)^\#
$$
defined on the higher-order Alexander module
$\SA:=H_1(M;\SR)$.
\end{thm}

\begin{proof} Note that $\SA$ is a torsion $\SR$--module by
Proposition~\ref{rank}, since $\SK$ is also the quotient
field of the Ore domain $\SR$. Define $B\l$ as the
composition of the following maps: the natural map
$\pi\co H_1(M;\SR)\lra H_1(M,\p M;\SR)$, the Poincar\'e duality
isomorphism to $H^2(M;\SR)$, the inverse of the Bockstein to
$H^1(M;\SK/\SR)$, and the usual Kronecker evaluation map to
$\SA^\#$. The Bockstein
$$
B\co H^1(M;\SK/\SR)\lra H^2(M;\SR)
$$
associated to the short exact sequence
$$
0\lra\SR\lra\SK\lra\SK/\SR\lra0
$$
is an isomorphism since $H^*(M;\SK)\cong H_*(M,\p
M;\SK)\cong0$ by Corollary~\ref{acyclic}.

We also need to show that $B\l$ is $``$conjugate symmetric".
The diagram below commutes up to a sign (see, for
example,~\cite[p.~410]{M}), where $B'$ is the homology
Bockstein
\begin{equation*}
\begin{split}
&\hskip100pt H_1(M;\SR)\\
&\hskip130pt \Bigg\downarrow\pi
\end{split}
\end{equation*}
\begin{equation}
\begin{split}
\begin{CD}
H_2(M,\p M; \SK/\SR)    @>B'>>  H_1(M,\p M;\SR)\\
@V\cong V{P.D.}V      @V\cong V{P.D.}V\\
H^1(M; \SK/\SR)    @>B>>    H^2(M;\SR)\\
@VV{\kappa}V\\
\end{CD}\\
\end{split}
\end{equation}
\hskip70pt$\Hom_{\SR}(H_1(M; \SR), \SK/\SR)$

\smallskip
\noindent and the two vertical homomorphisms are Poincar\'e
duality. Thus our map $B\l$ agrees with that obtained by going
counter-clockwise around the square and thus agrees with the
Blanchfield form defined by J.~Duval in a non commutative
setting~\cite[p.~623--624]{D}. The argument given there for
symmetry is written in sufficient generality to cover the
present situation and the reader is referred to it.
\end{proof}

\begin{defn}\label{Blanchfieldn} The {\it
$n^{\rm th}$--order linking form} for the knot $K$,
$B\l_n\co\SA_n(K)\to\SA_n(K)^\#$, is the linking form above
with $R=R_n$ (as in Section~\ref{rational}).
\end{defn}

\begin{prop}\label{non-singularity} The linking form
$B\l\co\SA_n(K)\lra\SA_n(K)^\#$ is non-singular if $n\neq1$. If
$n=1$ the kernel of $B\l$ is the submodule generated by the
longitude, and there is a non-singular linking form induced
on the ``reduced'' (quotient) module $\SA^*_1(K)$, obtained
by killing the longitude.
\end{prop}

\begin{cor}\label{duality} The localized modules $\SA_n(K)$
(if $n=1$ use $\SA^*_1(K)$) are self-dual. It follows that
the higher-order Alexander polynomials $e^n_i$ and $\Delta_n$
of Theorem~\ref{classification} are self-dual (an element
$\Delta$ of $R$ is {\it self-dual} if it is similar to
$\ov\Delta$).
\end{cor}

\begin{proof}[Proof of Corollary~\ref{duality}] Note that for
a finite cyclic module, $A=R/eR$, \\ $\Hom(A,\SK/R)\cong R/Re$ and
$A^\#\cong R/\bar e R$. The result then follows from the
uniqueness in Theorem~\ref{classification}.
\end{proof}

\begin{proof}[Proof of Proposition~\ref{non-singularity}] The
Kronecker map $H^1(\sk;R_n)\lra\linebreak
\Hom_{R_n}(\SA_n(K),\SK/R_n)$ is
an isomorphism since, over the PID $R_n$, the usual Universal
Coefficient Theorem holds (Remark~3.6.3) and
$\Ext_{R_n}(H_0(\sk;R_n),\linebreak
\SK_n/R_n)=0$ since $\SK/R_n$ is clearly
a divisible $R_n$--module and hence an injective $R_n$--module by
\cite[I Prop.~6.10]{Ste}. Thus $B\l$ is a isomorphism if and only
if the map $\pi\co H_1(\sk;R_n)\lra\ H_1(\sk,\p(\sk);R_n)$ is an
isomorphism. When $n=0$, the map $H_0(\p(\sk);\bbq[t,t^{-1}])
\lra\ H_0(\sk;\bbq[t,t^{-1}])$ is an isomorphism, implying $\pi$
is onto. Moreover $H_1(\p(\sk);\bbq[t^{\pm1}])$ has zero image in
$H_1(\sk;\bbq[t^{\pm1}])$ since any Seifert surface for $K$ lifts
to the $\G_0$ ($\infty$--cyclic) covering space, in other words the
longitude $\l\in\gs2$. Thus $\pi$ is an isomorphism when $n=0$.
Now suppose $n\ge2$. If $K$ has Alexander polynomial $1$ then
$\gs1$ is a perfect group so $\gs1=\gs n$ for all $n$ and thus
$\G_0=\G_n$ for all $n$ and $\SA_n=\SA_0$ for all $n$. The
non-singularity then follows from the $n=0$ case. Thus we may
assume that $\gs1/\gs2\neq0$. Below it will be shown that the
longitude is non-trivial in $\gs2/\gs3$. In particular the
longitude is non-trivial in $\G_n=G/\gns+$ if $n\ge2$. Since
$\l\in\gs1$, it follows that $\l$ is a non-trivial element of $\wt
G_{n+1}$. Therefore $\l-1$ is a non-zero element of $\bbz\wt
G_{n+1}$ and thus is invertible in $R_n$. Thus, by (the proof of)
Proposition~\ref{Ho}, $H_0(\p(\sk);R_n)=R_n/R_nI=0$, and so $\pi$
is surjective. Moreover since $\l$ is non-trivial in $\G_n$ if
$n\ge2$, $\pi_1(\p(\sk))$ embeds in $\G_n$ and the induced $\G_n$
cover is a union of planes so $H_1(\p(\sk);R_n)=0$ and $\pi$ is
also injective. This finishes the proof of the proposition in the
case $n\neq1$, modulo the proof that $\l\notin\gs3$ (assuming
$\gs1/\gs2\neq0$).

If $n=1$, the situation is more complicated. Let $E=\sk$ and
consider the commutative diagram below where all groups have
coefficients in $R_1$ unless specified.
\begin{equation*}
\begin{split}
&H_1(E)   \xrightarrow{\text{P.D.}}    H^2(E,\p E)
\xrightarrow{B^{-1}}
H^1(E,\p E;\SK_1/R_1)  \xrightarrow{\k}  \Hom(H_1(E,\p
E),\SK_1/R_1)\\
&\quad\Bigg\downarrow\pi_1 \hskip50pt
\Bigg\downarrow\pi_2\hskip70pt
\Bigg\downarrow\pi_3\hskip80pt
\Bigg\downarrow\pi^\#_1\\
&H_1(E, \p E)    \xrightarrow{\text{P.D.}}    H^2(E)
\xrightarrow{B^{-1}}    H^1(E;\SK_1/R_1)
\xrightarrow{\ \k\ }   \Hom(H_1(E),\SK_1/R_1)
\end{split}
\end{equation*}
All of the horizontal maps are isomorphisms. Let $f$ (respectively
$g$) denote the composition of all the maps in the top (bottom)
row. Then $B\l=g\circ\pi_1$ and its kernel is precisely the kernel
of $\pi_1$ which equals the image of $i_*\co H_1(\p E)\lra H_1(E)$.
This image is clearly generated by the longitude since the
infinite cyclic cover of $\p E$ is an annulus homotopy equivalent
to a circle representing a lift of the longitude. Moreover the
induced map $B\l^*$ is thus injective on $\SA^*_1\cong
H_1(E)/i_*(H_1(\p E))$. It remains only to show that the image of
$g\circ\pi_1$ is naturally isomorphic to $(\SA^*_1)^\#$, i.e.\ 
$\Hom(H_1(E)/H_1(\p E),\SK_1/R_1)$. The image of $g\circ\pi_1$
equals the image of $\pi^\#_1$ since $f$ is an isomorphism.
Consider the commutative diagram below. Since
$\ov\pi_1\co H_1(E)/H_1(\p E)\lra H_1(E,\p E)$
\begin{equation*}
\begin{split}
\Hom(H_1(E,\p E),\SK_1/R_1)\xrightarrow{\ov\pi^\#_1}
\Hom(H_1(E)/  &H_1(\p E),\SK_1/R_1)\\
\searrow\pi^\#_1\hskip70pt    &\Bigg\downarrow\pi^\#\\
\Hom(H_1  &(E),\SK_1/R_1)\\
&\Bigg\downarrow i^\#_*\\
\Hom(H_1  &(\p E),\SK_1/R_1)
\end{split}
\end{equation*}
is injective, its dual map $\ov\pi^\#_1$ (the horizontal map
above) is surjective since\linebreak
$\Ext_{R_1}(\un{\ \ },\SK_1/R_1)=0$ as remarked earlier in the
proof. Therefore the image of $\pi^\#$ is contained in the
image of $g\circ\pi_1$. Note that the image of (the diagonal
$\pi^\#_1$ is contained in the kernel of $i^\#_*$. But the
vertical sequence is exact and $\pi^\#$ is injective since
$\Hom$ is right exact. Thus
$\image(B^*\l)=\image(\pi^\#_1)=\image\pi^\#$ , and
$\pi^\#$ induces an isomorphism between $(\SA^*_1)^\#$ and
$\image(g\circ\pi_1)$. Therefore, with this identification,
$B\l$ induces an isomorphism $B\l^*$ between
$\SA^*_1\equiv\SA_1/\ker B\l$ and $(\SA^*_1)^\#$.
\end{proof}

\begin{prop}\label{longitude} If the (classical) Alexander
polynomial of $K$ is not $1$, then the longitude of $K$
represents a non-zero class in
$\gs2/\gs3\ox_{\bbz[G/\gs2]}R_1$. In particular
$\l\notin\gs3$.
\end{prop}

\begin{proof} Consider the coefficient system $\phi\co G\lra
G/\gs2\equiv\G_1$. Let $M$ be the result of zero framed
surgery on $K$ so $M=(\sk)\cup e^2\cup e^3$ where the
attaching circle of $e^2$ is the longitude. Since $\l\in\gs2$
for any knot, $\phi$ extends to $\pi_1(M)$. We may then
consider the commutative diagram of exact sequences below:
\begin{equation*}
\begin{split}
\begin{CD}
\bbk_1[t^{\pm1}]\\
@VV{\p'}V  \\
H_2((\sk)\cup e^2; \bbk_1[t^{\pm1}])   @>{\pi}>>
\bbk_1[t^{\pm1}]   @>{\p}>>  H_1(\sk;\bbk_1[t^{\pm1}])\\
@VV{i_*}V\\
H_2(M; \bbk_1[t^{\pm1}])  = 0\\
@VVV\\
0\\
\end{CD}\\
\end{split}
\end{equation*}
The horizontal sequence is that of the pair $(\sk\cup
e^2,\sk)$ and the generator of the $\bbk_1[t^{\pm1}]$ may be
thought of as $e^2$ and its boundary as the class represented
by the longitude in $\SA_1$. Suppose
$\l\in\SA_1\equiv\gs2/\gs3\ox_{\bbz[G/\gs2]}R_1$ is zero.
Then the map $\pi$ would be a surjection. Now consider the
vertical exact sequence of the pair $(M,\sk\cup e^2)$. Here
the generator of $\bbk_1[t^{\pm1}]$ may be thought of as the
$3$--cell $e^3$. We have $H_2(M;\bbk_1[t^{\pm1}])\cong H^1(M;
\bbk_1[t^{\pm1}])\cong\Ext(H_0(M;\bbk_1[t^{\pm1}]),\bbk_1[t^{\pm1}])$.
If the Alexander polynomial of $K$ is not $1$ then the
Alexander module $\gs1/\gs2$ contains some $x\neq e$. Thus
$x-1$ lies in the augmentation ideal of $\bbz G$ and
$\phi(x-1)$ is invertible in $\bbk_1[t^{\pm1}]\equiv R_1$
since $x\in\wt G$ (see Proposition~\ref{gamma-n}). Thus
$H_0(M;\bbk_1[t^{\pm1}])$ vanishes by (the proof of)
Proposition~\ref{Ho} and hence $H_2(M;\bbk_1[t^{\pm1}])=0$.
Therefore $\p'$ and $\pi\circ\p'$ are epimorphisms. We claim
that the diagonal map $(\pi\circ\p')$ sends $1\to1-t$. This
claim is seen by analyzing how the $3$--cell goes over the
$2$--cell twice. This map is clearly not surjective since
$1-t$ is not a unit. This contradicts our assumption that the
longitude vanished.
\end{proof}

\section{Calculation from a presentation of the knot
group}\label{fox}

A presentation matrix for $\SA_n(K)$ can be derived from any
finite presentation of $G=\pi_1(\sk)$.

It is known that, for any regular $\G$ covering space $X_\G\to X$
of a finite complex, the free differential calculus can be used to
give a presentation matrix for $H_1(X_\G,\tl x_0)$ as a
$\BZ\G$--module where $\tl x_0$ is the inverse image of a basepoint
(see, for example \cite{H}). The torsion submodule of
$H_1(X_\G,\tl x_0)$ can easily seem to be isomorphic to
$H_1(X_\G)$. Thus a presentation matrix can be computed for a
module whose torsion submodule is $\SA^\bbz_n(K)$. The same holds
for $\SA_n(K)$. Over a PID, it is theoretically possible to
simplify a presentation matrix by appropriate row and column
operations until it is diagonal, thus calculating the $\d_n$
(see\cite{Ha}). This necessitates deciding whether or not a given
element of the solvable group $\gs n/\gns+$ is trivial. Sometimes
this is difficult. However note that for $n=1$ this quotient group
is merely the classical Alexander module of the knot. Hence there
exists a practical algorithm to compute $\SA_1(K)$. We hope to
soon implement this. Details and some sample calculations can be
found in \cite{Ha}.

\section{Questions and open directions}\label{Questions}

\begin{enumerate}
\item Find invariants of the higher-order modules which can
detect the peripheral structure of a knot.
\item Find other invariants of the integral modules that are
not simply invariants of the localized modules.
\item Develop effective invariants of the higher-order
Alexander polynomials or find ways to reduce their
indeterminacy.
\item Is there a higher-order Seifert form? (The existence of
$(t-1)$--torsion has thwarted our efforts on this question.)
\item Is there a knot $K$ and some $n>0$ for which $\d_n(K)$
is a non-zero {\it even} integer? If not then a complete
realization theorem for the $\d_i$ can be derived from the
techniques of Section~\ref{examples}.
\item Find higher-order
Seiberg-Witten invariants of $3$--manifolds that reflect these
higher-order modules.
\item Are the invariants $\d_i$ of finite type?
\item Prove that for each $n\ge0$ there exist
knots $K$ and $K_*$ such that $\SA_i(K)\cong\SA_i(K_*)$ for
$0\le i<n$ but $\b_n(K_*)\neq\b_n(K)$.
\item The Arf invariant of a knot is determined by its
Alexander polynomial which is in turn determined by its
Alexander module which is in turn determined by any Seifert
matrix. Similarly the Levine-Tristram signatures of a knot
are determined by the Alexander module and its Blanchfield
form which are in turn determined by a Seifert matrix. Can
any such statements be made for the higher-order bordism
invariants $\b_n$, modules $\SA^\bbz_n$, signatures $\rho_n$
and presentation matrices from Section~\ref{seifert}?
\item Find knots with the same higher-order modules but
different linking forms.
\item Find ways to compute the $\rho_n$.
\item Apply these ideas to links, string links, braids and
mapping class groups.
\item Do these invariants have any special behavior on other
special classes of knots? (for example connected-sums of
knots have non-longitudinal $(t-1)$--torsion in $\SA_1$).
\item Find applications to contact structures on $3$--manifolds
(which seem to be closely related to fibering questions).
\end{enumerate}

\Addresses\recd

\begin{thebibliography}{LMO}

\bibitem[A]{A} M.F. Atiyah, {\it Elliptic operators, discrete
groups and von Neumann algebras}, Asterisque 32--33, 43--72,
1976.

\bibitem[BW]{BW} M. Boileau and S. Wang, {\it Non-zero degree
maps and surface bundles over $S^1$}, J. Diff. Geometry 43,
789--806, 1996.

\bibitem[Br]{Br} K.S. Brown, {\it Cohomology of Groups}
Springer-Verlag 1982.

\bibitem[Ch]{Ch} Jae Choon Cha, {\it Fibred Knots and Twisted
Alexander Invariants}, preprint 2001, \arxiv{math.GT/0109136}.

\bibitem[CO]{CO} T. Cochran and K. Orr, {\it Homology boundary links and
Blanchfield forms: Concordance classification and new
tangle-theoretic constructions}, Topology 33, no. 3, 397--427,
1994.

\bibitem[COT1]{COT1} T. Cochran, K. Orr and P. Teichner,
{\it Knot concordance, Whitney towers and $L^2$--signatures},
Annals of Math.,157 ,433--519, 2003.

\bibitem[COT2]{COT2} T. Cochran, K. Orr and P. Teichner, {\it
Structure in the classical knot concordance group}, Commentarii
Math. Helvetici, 79, 105--123, 2004.

\bibitem[CT]{CT} T. Cochran  and P. Teichner, {\it Knot concordance and von
Neumann eta invariants}, preprint, 2002, Rice University.

\bibitem[Co1]{Co1} P.M. Cohn,  {\it Skew Fields}, Cambridge Univ. Press,
Cambridge 1995.

\bibitem[Co2]{Co2} P. M. Cohn,  {\it Free Rings and their Relations}, Second
Edition, London Math. Soc. Monographs no. 19, Academic Press, London and
New York 1985.

\bibitem[ChG]{ChG} J. Cheeger and M. Gromov, {\it Bounds on
the von Neumann dimension of $L^2$--cohomology and the Gauss-Bonnet
Theorem for open manifolds}, J. Diff. Geometry 21, 1--34, 1985.

\bibitem[CTe]{CTe} J. Conant and P. Teichner, {\it Grope
cobordism of classical knots}, Topology 43, no. 1, 119--156, 2004.

\bibitem[Do]{Do}  Dodziu, Linnell, Matai, Schick, Yates, {\it Approximating
L2-invariants, and the Atiyah Conjecture},Dedicated to the memory
of J\"orgen K. Moser.  Comm. Pure Appl. Math.  56,  no. 7, 839--873,
2003.

\bibitem[D]{D}  J. Duval, {\it Forme de Blanchfield et cobordisme
d'entrelacs bords}, Commentarii Math. Helvetici 61, 617--635,
1986.

\bibitem[GR]{GR} S. Garoufalidis and L. Rozansky, {\it The loop expansion of
the Kontsevich integral, abelian invariants of knots and
$s$--equivalence}, \arxiv{math.GT/0003187}, Topology (in
press).

\bibitem[GL]{GL} S. Garoufalidis and J. Levine, {\it Homology surgery and
invariants of $3$--manifolds}, \href{http://www.maths.warwick.ac.uk/gt/GTVol5/paper18.abs.html}{Geometry and Topology, vol 5,
551--578, 2001.}

\bibitem[Go]{Go} C. Gordon, {\it Some aspects of classical knot theory}  Knot theory (Proc. Sem., Plans-sur-Bex, 1977),
  pp. 1--60, Lecture Notes in Math., 685, Springer, Berlin, 1978.

\bibitem[Hb]{Hb} K. Habiro, {\it Claspers and finite type invariants of
links}, Geometry and Topology, \href{http://www.maths.warwick.ac.uk/gt/GTVol4/paper1.abs.html}{vol 4, 1--83, 2001.}

\bibitem[Ha]{Ha} S. Harvey, {\it Higher Order Polynomial
Invariants of $3$--manifolds giving lower bounds for the Thurston
Norm}, preprint, Rice University 2001, \arxiv{math.GT/0207014}.

\bibitem[H1]{H1}  J. Harer,  {\it How to construct all fibered
knots and links}, Topology 21 no.~3, 263--280, 1982.

\bibitem[H]{H} J. Hempel, {\it Intersection Calculus on
Surfaces with Applications to $3$--manifolds}, Memoirs of AMS
\#282, volume~43, American Mathematical Society, Providence,
Rhode Island (1983).

\bibitem[Hig]{Hig} G. Higman, {\it The units of group rings},
Proc. London Math. Soc. 46, 231--248, 1940.

\bibitem[Hi2]{Hi2} J. Hillman, {\it Alexander Ideals of
Links}, LNM \#895, Springer-Verlag, Berlin (1981).

\bibitem[HK]{HK} N. Higson and Kasparov, Preprint 1998.

\bibitem[HS]{HS} P.J. Hilton and U. Stammbach, {\it A Course
in Homological Algebra} Springer-Verlag, 1971, New York.

\bibitem[HoS]{HoS} J. Howie and H.R. Scheebeli, {\it
Homological and topological properties of locally indicable
groups}, Manuscripta Math. 44, no. 1--3, 71--93, 1983.

\bibitem[J]{J} N. Jacobson, {\it The Theory of Rings},
American Math. Soc., New York, 1943.

\bibitem[Ka]{Ka} Akio Kawauchi, {\it Almost identical
imitations of $(3,1)$--dimensional manifold pairs}, Osaka J. Math.
26, 743--758, 1989.

\bibitem[Ki]{Ki} R. Kirby, {\it The calculus of framed links
in $S^3$}, Inventiones Math. l45, 35--56, 1978.

\bibitem[KL]{KL} P. Kirk and C. Livingston, {\it Twisted
Alexander invariants, Reidemeister torsion and Casson-Gordon
invariants}, Topology 38, no. 3, 635--661, 1999.

\bibitem[Kr]{Kr} P. Kronheimer, {\it Embedded Surfaces and
Gauge Theory in Three and Four Dimensions}, Surveys in
differential geometry, Vol. III, Int. Press, 243--298, 1998.

\bibitem[Kr2]{Kr2} P. Kronheimer, {\it Minimal genus in
$S^1\x M^3$}, Inventiones Math. l34, 363--400, 1998.



\bibitem[L1]{Le} J. P. Levine, {\it Knot cobordism groups in
codimension two}, Comm. Math. Helv. 44, 229--244, 1969.

\bibitem[L3]{L3} J. P. Levine, {\it Knot Modules, I}, Trans.
Amer. Math. Soc. 229, 1977, 1--50.

\bibitem[Lew]{Lew} J. Lewin, {\it A note on zero divisors in
group rings}, Proc. Amer. Math. Soc. 31, 1972, 357--359.

\bibitem[Ma]{Ma} V. Mathai, {\it Spectral flow, eta invariants and von
Neumann algebras}, J. Functional Analysis 109, 442--456, 1992.

\bibitem[M]{M} J. Munkres, {\it Elements of Algebraic Topology},
Addison-Wesley, Menlo Park CA, 1984.

\bibitem[M2]{M2} J. Munkres, {\it Topology}, 2nd edition,
Prentice Hall, Saddle River, N.J. (2000).

\bibitem[P]{P} D. Passman, {\it The Algebraic Structure of
Group Rings}, John  Wiley and Sons, 1977, New York.

\bibitem[NS]{NS} S. Naik and T. Stanford, {\it A move on
diagrams that generates $S$--equivalence of knots},
\arxiv{math.GT/9911005}, to appear in J. Knot Theory and its
Applications.

\bibitem[Ra2]{Ra2} A. Ranicki, {\it Exact Sequences in the
Algebraic Theory of Surgery}, Princeton University Press,
1981, Princeton New Jersey.


\bibitem[Ru]{Ru} S. Roushon, {\it Topology of $3$--manifolds and a class of
groups}, preprint 2002, \arxiv{math.GT/0209121}.


\bibitem[Sta]{Sta} J. Stallings, {\it Constructions of fibered knots and
links}, AMS Proceedings of Symposia in Pure Math. 32, 1978, 55--60.

\bibitem[St]{St} J. Stallings, {\it Homology and Central series of Groups},

\bibitem[Ste]{Ste} B. Stenstr\"om, {\it Rings of Quotients},
Springer-Verlag, 1975, New York.

\bibitem[Sti]{Sti} A. Stoimenow, {\it Vassiliev invariants on
fibered and mutually obverse knots}, J. Knot Thry. Ram. 8
(1999), no. 4, 511--519.

\bibitem[Str]{Str} R. Strebel, {\it Homological methods applied to the
derived series of groups}, Comment. Math. 49, 1974, 302--332.

\bibitem[W]{W} G.W. Whitehead, Elements of Homotopy Theory,
{\it Graduate Texts in Mathematics \#61}, Springer-Verlag, 1978, New York.

\bibitem[Wa]{Wa} C. T. C. Wall, {\it Surgery on Compact
Manifolds}, London Math. Soc. Monographs~1, Academic Press~1970.


\end{thebibliography}
\end{document}